\documentclass[final,hidelinks,onefignum,onetabnum]{siamart220329}

\usepackage{lipsum}
\usepackage{amsfonts}
\usepackage{graphicx}
\usepackage{epstopdf}
\usepackage{algorithmic}
\ifpdf
  \DeclareGraphicsExtensions{.eps,.pdf,.png,.jpg}
\else
  \DeclareGraphicsExtensions{.eps}
\fi

\newsiamremark{remark}{Remark}
\newsiamremark{hypothesis}{Hypothesis}
\crefname{hypothesis}{Hypothesis}{Hypotheses}
\newsiamthm{claim}{Claim}

\headers{Solving Kolmogorov backward equations with FHT operators}{X. Tang, L. Collis, L. Ying}

\title{Solving high-dimensional Kolmogorov backward equations with functional hierarchical tensor operators }

\author{Xun Tang\thanks{Institute for Computational and Mathematical Engineering, Stanford, CA 94305, USA. 
  (\email{xuntang@stanford.edu}).}
  \and
  Leah Collis\thanks{Institute for Computational and Mathematical Engineering, Stanford, CA 94305, USA. 
  (\email{lcollis@stanford.edu}).}
  \and
  Lexing Ying\thanks{Department of Mathematics and Institute for Computational and Mathematical Engineering, Stanford, CA 94305, USA. 
  (\email{lexing@stanford.edu})}
  \funding{This work is partially supported by NSF grants DMS-2011699 and DMS-2208163.}
  }

\usepackage{amsopn}

\usepackage{amsmath,amssymb,a4wide}
\usepackage{latexsym}
\usepackage{indentfirst}
\usepackage{graphicx}
\usepackage{placeins}
\usepackage{enumerate}
\usepackage{booktabs}
\usepackage{algorithm}
\usepackage{algorithmic}
\usepackage{multirow}
\usepackage{optidef}
\usepackage{graphicx,color}
\usepackage{diagbox}
\usepackage[normalem]{ulem}

\ifpdf
\hypersetup{
  pdftitle={KBE-with-FHT},
  pdfauthor={X. Tang, L. Collis, and L. Ying}
}
\fi

\newcommand{\R}{\mathbb{R}}

\begin{document}

\maketitle

\begin{abstract}
Solving high-dimensional partial differential equations necessitates methods free of exponential scaling in the dimension of the problem. This work introduces a tensor network approach for the Kolmogorov backward equation via approximating directly the Markov operator. We show that the high-dimensional Markov operator can be obtained under a functional hierarchical tensor (FHT) ansatz with a hierarchical sketching algorithm. When the terminal condition admits an FHT ansatz, the proposed operator outputs an FHT ansatz for the PDE solution through an efficient functional tensor network contraction procedure. In addition, the proposed operator-based approach also provides an efficient way to solve the Kolmogorov forward equation when the initial distribution is in an FHT ansatz. We apply the proposed approach successfully to two challenging time-dependent Ginzburg-Landau models with hundreds of variables. 
\end{abstract}




\section{Introduction}

This paper focuses on solving the Kolmogorov backward equation (KBE):
\begin{equation}\label{eqn: KBE}
    \begin{cases}
        &\partial_{t}u -  \nabla V \cdot  \nabla u + \frac{1}{\beta}\Delta u = 0,
    \quad x \in \Omega \subset \R^{d}, t \in [0, T),\\
    &u(y, T) = f(y), \quad  y \in \Omega \subset \R^{d},
    \end{cases}
\end{equation}
where $V$ is a potential function and $\beta$ is the inverse temperature. The solution \(u\) can also be represented as \(u(x, t) = \mathbb{E}\left[f(X_{T}) \mid X_{t} = x\right] \), where 
$X_t$ follows the Langevin dynamics 
\begin{equation}\label{eqn: langevin}
    dX_t = -\nabla V(X_t) dt + \sqrt{2\beta^{-1}}dB_{t}.
\end{equation}

In this work, we focus on the high-dimensional setting where the domain has hundreds of dimensions. One prevalent type of high-dimensional Kolmogorov backward equation comes from the discretization of an infinite-dimensional Markov chain into finite but very large dimensions. Similar to the prior work in \cite{tang2023solving}, we focus on the motivating example of a time-dependent Ginzburg-Landau (G-L) model used for studying the phenomenological theory of superconductivity \cite{ginzburg2009theory, hoffmann2012ginzburg, hohenberg2015introduction,weinan2004minimum}. 
The 2D G-L model considers a particle taking the form of a field \(x(a) \colon [0,1]^2 \to \R\). The potential function is a functional $V(x)$ defined as follows:
\begin{equation}\label{eqn: GL potential infinite}
  V(x) = \frac{\lambda}{2}\int_{[0,1]^2} |\nabla_a x(a)|^2 \, da + \frac{1}{4\lambda}\int_{[0,1]^2} |1- x(a)^2|^2 \, da,
\end{equation}
where \(\lambda\) is a parameter balancing the relative importance of the first and second term in \eqref{eqn: GL potential infinite}.
Under discretization, one decomposes the unit square \([0,1]^2\) into a grid of \(d = m^2\) points $\{(ih,jh)\}$, for
$h=\frac{1}{m+1}$ and $1\le i,j\le m$. The discretized field is the $d$-dimensional vector $x=(x_{(i,j)})_{1\le i,j \le m}$, where \(x_{(i, j)} = x(ih, jh)\). The associated discretized potential energy is defined as
\begin{equation}\label{eqn: 2D GZ model}
V(x) = V(x_{(1,1)}, \ldots ,x_{(m,m)}) := \frac{\lambda}{2} \sum_{v \sim w}\left(\frac{x_{v} - x_{w}}{h}\right)^2 +  \frac{1}{4\lambda} \sum_{v}\left(1 - x_v^2\right)^2,
\end{equation}
where $v$ and $w$ are Cartesian grid points and \(v \sim w\) if and only if they are adjacent. For any \(f \colon \mathbb{R}^{d} \to \mathbb{R}\), the evolution of the conditional expectation with \(f\) as a terminal condition constitutes the Kolmogorov backward equation \eqref{eqn: KBE} with the discretized potential function $V(x)$. Even a moderate discretization scheme leads to a high-dimensional problem.

\subsection{Background and contribution}\label{sec: contribution}\paragraph{Existing methods for solving the Kolmogorov backward equation} 
While the traditional finite difference \cite{brennan1978finite, zhao2007compact} and finite element schemes \cite{kumar2006solution} are candidate numerical solvers for the Kolmogorov backward equation, they cannot be used for the high-dimensional case due to the curse of dimensionality. More recent work \cite{berner2020numerically, beck2021solving} shows that one can overcome the curse of dimensionality by combining Monte-Carlo methods with a downstream regression task. Essentially, due to the conditional expectation formula, one can use Monte-Carlo to approximately query function evaluation of \(u\) at a sufficiently large number of points in \(\Omega \times [0, T]\). After obtaining the training sample, one can use the approximate value of \(u\) at these points as training data for a neural network regression task. 

While the regression-based perspective overcomes the curse of dimensionality, its implementation faces two major challenges. First, the method incurs a significant cost to prepare the training target, as a large number of SDE simulations is required for each training sample. Second, retraining of the neural network is unavoidable when one considers a different terminal condition for \(u(\cdot, T)\). 

\paragraph{Markov operator for solving the Kolmogorov backward equation}
In this work, we propose to solve the Kolmogorov backward equation by approximating the solution operator \cite{lin2011fast,fan2020solving,kovachki2021neural}, i.e., obtaining a mapping from the terminal condition \(f\) to the solution \(u\). For any \(t \in [0, T]\) and a stochastic process \(\left(X_{s}\right)_{s \in [0, T]}\) satisfying \eqref{eqn: langevin}, one has the following equation:
\begin{equation}\label{eqn: Bayes}
    u(x, T-t) = \int_{\Omega} \mathbb{P}\left[X_T = y \mid X_{T-t} = x \right]f(y) \, dy =\int_{\Omega} \mathbb{P}\left[X_{t} = y \mid X_{0} = x \right]f(y) \, dy
\end{equation}
where the second equality of \eqref{eqn: Bayes} is due to the time invariance of the dynamics for \eqref{eqn: langevin}.
Thus, the solution operator \(G_t \colon \Omega \times \Omega \to \mathbb{R}\) is simply the conditional density defined by 
\begin{equation}\label{eqn: Markov operator}
    G_t(x, y) := \mathbb{P}\left[X_{t} = y \mid X_0 = x \right],
\end{equation}
which is typically referred to as the \emph{Markov operator}.

Our proposed approach to approximate the Markov operator is motivated by Bayes' rule. For any chosen initial distribution \(X_0 \sim p_0\), one has 
\begin{equation}\label{eqn: Markov operator by Bayes}
    G_{t}(x, y) = \frac{1}{p_0(x)}\mathbb{P}\left[X_{0} = x, X_{t} = y \right],
\end{equation}
which shows that the conditional density \(G_{t}\) can be characterized through the probability density of the joint variable \((X_0, X_{t})\), which is an important observation underlying our subsequent tensor network approach to represent the Markov operator.




\paragraph{Functional hierarchical tensor operator} 
Functional hierarchical tensor (FHT) \cite{tang2023solving,hackbusch2012tensor,hackbusch2009new} is a class of parametric functions suitable for reduced-order modeling. While a detailed introduction can be found in \cite{tang2023solving}, we mention a few of its appealing properties. The FHT ansatz scales linearly in dimension for memory requirement, function evaluation, moment calculation, and sampling (when used as an ansatz for density). Most importantly, the FHT ansatz scales linearly in dimension to calculate marginal distribution density and conditional distribution density, which sets it apart from common neural network architectures for density estimation.

Due to the aforementioned property, one can use an FHT ansatz to represent the probability density of the joint variable \((X_{0}, X_{t})\), which leads to an FHT-based approximation to the operator \(G_{t}\).
The learned operator allows an efficient way to obtain the approximate value of \(u(x, T-t)\) for arbitrary point \(x \in \Omega\) and terminal condition \(f\). When \(f\) has an FHT structure, equation \eqref{eqn: Bayes} allows for more efficient computation of \(u(x, T-t)\) through a tensor diagram contraction.

\paragraph{Density estimation through hierarchical sketching}
Obtaining the joint density distribution is done by using a modification of an FHT-based approach in solving the Fokker-Planck equations \cite{tang2023solving}. After running particle simulations over the Langevin dynamics in equation \eqref{eqn: langevin}, one collects the value of the particles at the initial time \(0\) and the intermediate time \(t\). Then, one performs a hierarchical sketch-based density estimation subroutine \cite{peng2023generative,tang2023solving} on samples of the joint variable \((X_{0}, X_{t})\). The outputted joint density is used to represent the Markov operator according to \eqref{eqn: Markov operator by Bayes}. 

Compared to the approach in \cite{tang2023solving}, the hierarchical sketching subroutine in this instance has new challenges, such as the design of the hierarchical bipartition structure and the functional basis. Using sophisticated numerical techniques, we are able to solve the Ginzburg-Landau model with \(d\) in the hundreds.

\subsection{Related work}\label{sec: related}

\paragraph{Neural network methods for density estimation}
As the Markov operator is obtained through density estimation, our work is related to generative learning algorithms for which density estimation is possible. Representing \(G_{t}\) with the joint distribution function for \((X_{0}, X_{t})\) is related to density estimation. Examples include normalizing flows \cite{tabak2010density,rezende2015variational}, energy-based models \cite{hinton2002training, lecun2006tutorial}, and the more recent diffusion models \cite{song2019generative,song2021maximum}. The main challenge in obtaining the Markov operator through the aforementioned works lies in the calculation of the normalization constant and high-dimensional integration through \eqref{eqn: Bayes}.
On the other hand, representing \(G_{t}\) directly as a conditional distribution is related to probabilistic imputation methods \cite{yoon2018estimating, cao2018brits,che2018recurrent} and conditional generative modeling \cite{mirza2014conditional,sohn2015learning,winkler2019learning,bhattacharyya2019conditional,lu2020structured}. For neural network-based density estimation tasks, the training associated with a neural network might be slow, whereas our sketch-based approach relies on simple linear algebraic subroutines and does not involve an iterative training component.

\paragraph{Monte-Carlo methods} 
As one has access to trajectories of the particle dynamics due to Monte-Carlo sampling, it is possible to perform sampling to obtain the value of \(u\) at arbitrary points. Our proposed sketching approach approximates the Markov operator through data assimilation of one single batch of particle simulations, while the Monte-Carlo approach requires a new set of particle simulations for every new point queried for the PDE solution \(u\). Thus, the Monte-Carlo approach is less suitable when one is interested in the solution at more than a few points or when one is interested in the global properties of the solution.

\paragraph{Tensor network methods}
We refer the readers to \cite{tang2023solving} for a detailed discussion on tensor network-based approaches for solving high-dimensional PDEs. It is worth mentioning that a sketch-based density estimation subroutine for the functional tensor train (FTT) has been developed in \cite{hur2023generative}. The proposed approach of this work can be similarly adapted to approximate the Markov operator with a functional tensor train as the tensor network architecture. However, the FTT ansatz is less capable of capturing the global correlation structure, which is why this work instead focuses on the FHT ansatz. It is also possible to use the regression-based perspective mentioned in Section \ref{sec: contribution} with the functional tensor network as the parameterized function class for performing the regression task, and a similar approach has been taken in \cite{richter2021solving, chen2023committor}. 

\subsection{Contents and notations} \label{sec: notation}

We outline the structure of the remainder of the manuscript.
Section \ref{sec: background} details the functional hierarchical tensor structure for approximating the Markov operator.
Section \ref{sec: alg} goes through the algorithm for solving the Kolmogorov backward equation with the functional hierarchical tensor. 
Section \ref{sec: numerics} details the numerical implementations and results from solving Ginzburg-Landau models in 1D and 2D. 

For notational compactness, we introduce several shorthand notations for simple derivation. For \(n \in \mathbb{N}\), let \([n] := \{1,\ldots, n\}\). For an index set \(S \subset [d]\), we let \(x_{S}\) stand for the subvector with entries from index set $S$. We use \(\Bar{S}\) to denote the set-theoretic complement of $S$, i.e. \(\Bar{S} = [d] - S\). 


\section{Tensor network structure for the Markov operator}\label{sec: background}
To address the Ginzburg-Landau and general models coming from discretizing infinite-dimensional Markov chains, we propose to represent the Markov operator via a functional hierarchical tensor ansatz based on a binary-tree-based low-rank structure. 
In this section, we go over the functional hierarchical tensor and the tensor network structure for the Markov operator. Importantly, to simplify the discussion, we discuss obtaining the Markov operator for a fixed \(t \in [0, T]\).
The symbol \(d\) is reserved for the dimension of the state-space variable, and the symbol \(N\) is reserved for the number of samples. 

\subsection{Functional hierarchical tensor}

For completeness, we give a short introduction to the functional hierarchical tensor ansatz, which largely follows the more detailed exposition in \cite{tang2023solving}. We first describe the \emph{functional tensor network} representation of a \(2d\)-dimensional function \(g \colon \mathbb{R}^{2d} \to \mathbb{R}\). 
Let \(\{\psi_{i;j}\}_{i = 1}^{n}\) denote a collection of orthonormal function basis over a single variable for the \(j\)-th co-ordinate, and let \(C \in \R^{n^{2d}}\) be the tensor represented by a tensor network. The functional tensor network is the \(2d\)-dimensional function defined by the following equation:
\begin{equation}\label{eqn: htn forward map}
    g(x)\equiv g(z_{1}, \ldots, z_{2d}) = \sum_{i_{1}, \ldots, i_{2d} = 0}^{n-1} C_{i_1,\ldots, i_{2d}} \psi_{i_1; 1}(z_1)\cdots \psi_{i_{2d}; 2d}(z_{2d}) = \left<C, \, \bigotimes_{j=1}^{2d} \Vec{\Psi}_{j}(z_j) \right>,
\end{equation}
where \(\Vec{\Psi}_{j}(z_j) = \left[\psi_{1;j}(z_j),\ldots, \psi_{n;j}(z_j) \right]\) is an \(n\)-vector encoding the evaluation for the \(z_j\) variable over the entire \(j\)-th functional basis.
In particular, \(g\) is a \emph{functional hierarchical tensor} when \(C\) is represented by a hierarchical tensor network ansatz.

\begin{figure}[!ht]
    \centering
    \begin{minipage}{0.48\textwidth}
        \centering
        \includegraphics[width=\textwidth]{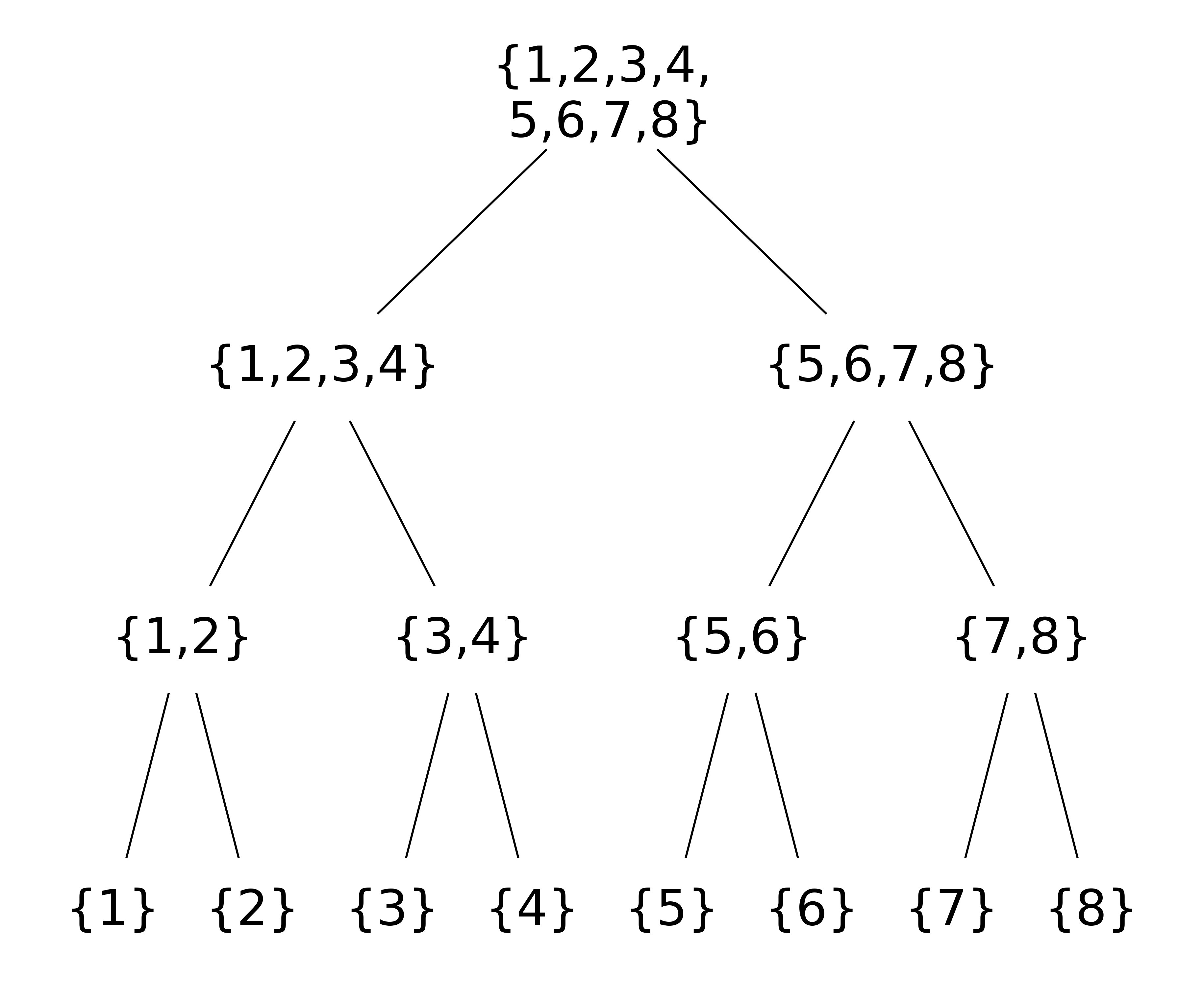}
        \caption{Binary decomposition of variables}
        \label{fig:binary_tree_8_nodes_subfig}
    \end{minipage}\hfill
    \begin{minipage}{0.50\textwidth}
        \centering
        \includegraphics[width=\textwidth]{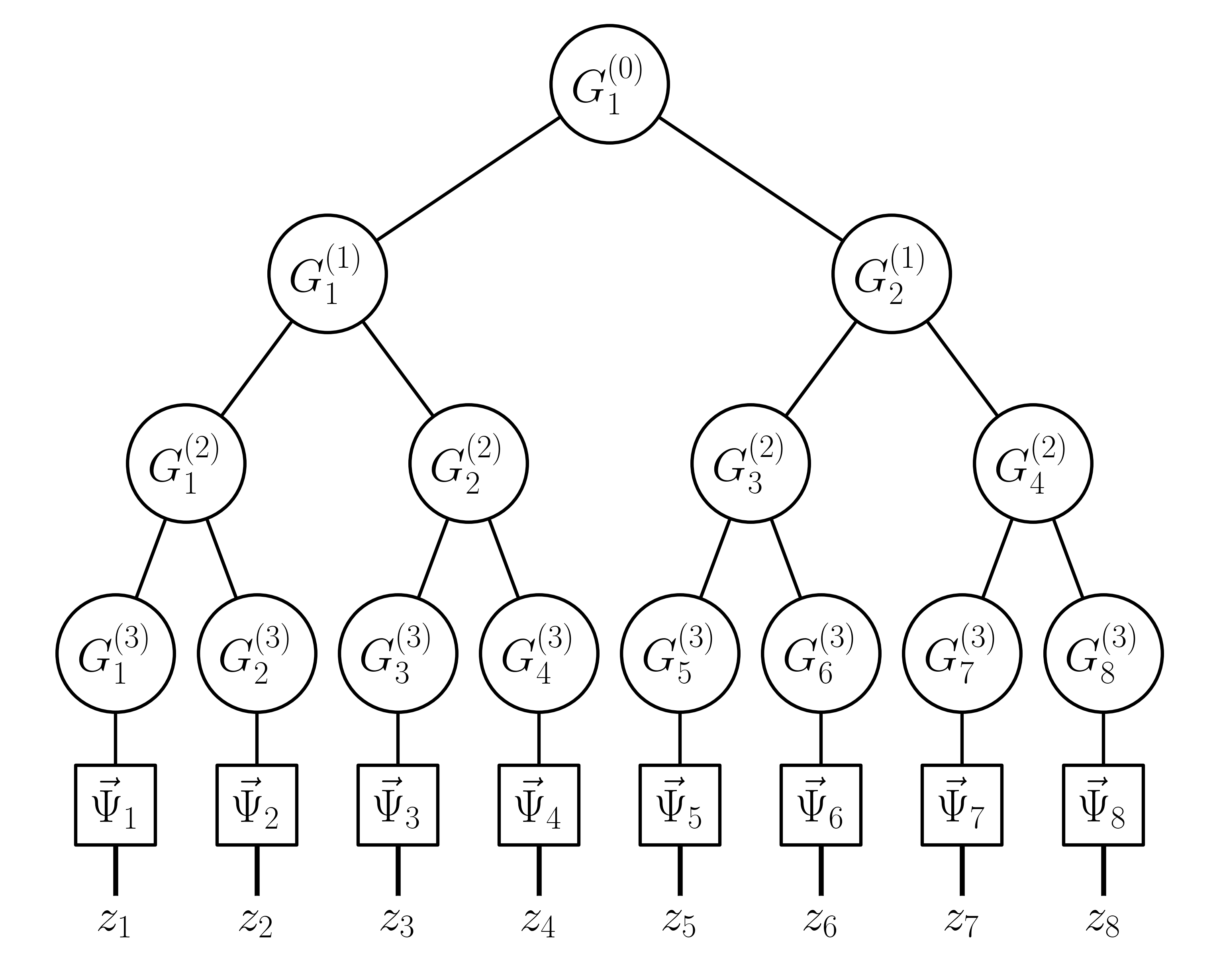}
        \caption{Functional hierarchical tensor diagram}
        \label{fig:FHT_L_3}
    \end{minipage}
    \caption{Illustrations of functional hierarchical tensor with $d = 4$.}
    \label{fig:decompositions}
\end{figure}

The hierarchical tensor network is characterized by a hierarchical bipartition of the variable set. Without loss of generality, let \(d = 2^{L}\) so that the variable set admits exactly \(L+1\) levels of variable bipartition. As illustrated in Figure \ref{fig:binary_tree_8_nodes_subfig}, at the \(l\)-th level, the variable index set is partitioned according to 
\begin{equation}\label{eqn: bipartition}
    [2d] = \bigcup_{k = 1}^{2^{l}} I_{k}^{(l)}, \quad I_{k}^{(l)} := \{ 2^{L - l + 1}(k-1) + 1, \ldots, 2^{L - l + 1}k\},
\end{equation}
which in particular implies the recursive relation that \(I_{k}^{(l)} = I_{2k-1}^{(l+1)}\cup I_{2k}^{(l+1)}\). 

Given the hierarchical bipartition, the tensor network structure of a hierarchical tensor network has the same tree structure as in the variable bipartition, which is illustrated by Figure \ref{fig:decompositions}.
Moreover, as Figure \ref{fig:decompositions} shows, each internal bond of 
the tensor network is associated with an edge in the binary decomposition of the variables. Each physical bond of the tensor network connects a leaf node with a $\Vec{\Psi}_j(z_j)$ node. In summary, the graphical structure of the hierarchical tensor network is formed by a binary tree, and the physical index is at the leaf node. In terms of runtime complexity and memory complexity, the functional hierarchical tensor enjoys the same \(O(d)\) scaling as that of a hierarchical tensor \cite{peng2023generative}. 

\subsection{Main idea for the FHT representation of the Markov operator}

We now describe the functional hierarchical tensor representation of the Markov operator. 
Let the initial distribution for \(\left(X_{s}\right)_{s \in [0, T]}\) satisfy \(X_0 \sim p_0\). By Bayes' rule and time invaraince of \eqref{eqn: langevin}, the Markov operator \(G_t = \mathbb{P}\left[X_{t} = y \mid X_{0} = x \right]\) satisfies the equation 
\[
    G_{t}(x, y) := \frac{1}{p_0(x)}\mathbb{P}_{X_0 \sim p_0}\left[X_{0} = x, X_{t} = y \right],
\]
where \(X_0 \sim p_0\) means the stochastic process \(\left(X_{s}\right)_{s \in [0, T]}\) follows from the Langevin dynamics in \eqref{eqn: langevin} with initial condition \(X_0 \sim p_0\).

In practice, one solves the KBE solution \(u\) inside a bounded domain \(\mathcal{X} \subset \Omega\), and correspondingly one considers the value of \(G_{t}(x, y)\) inside \(x \in \mathcal{X}\). Thus, one can set \(p_0\) as the uniform distribution on \(\mathcal{X}\). Under this construction, one has
\begin{equation}\label{eqn: Markov operator new formula}
    G_{t}(x, y) = \mathrm{Vol}(\mathcal{X})\mathbb{P}\left[X_{0} = x, X_{t} = y \right], \quad x \in \mathcal{X}.
\end{equation}

In particular, the construction of \(G_{t}\) in \eqref{eqn: Markov operator new formula} implies solving the Markov operator can be reduced to density estimation over the joint variable \((X_0, X_{t})\). Thus, one can employ the sketch-based density estimation subroutine in \cite{tang2023solving} to estimate the probability distribution of \((X_0, X_{t})\) as a functional hierarchical tensor ansatz, which gives the Markov operator an FHT representation. The rest of this section gives implementation details to this approach.

\subsection{Network structure of the Markov operator}
Representation of the Markov operator in the FHT ansatz leads to the question of the hierarchical tensor network structure for the \(2d\)-dimensional probability density function
\[g(x, y) = \mathbb{P}_{X_0 \sim p_0}\left[X_{0} = x, X_{t} = y \right].\]

We propose to use the interlacing scheme to order the joint variable. For \(x = (x_1, \ldots, x_d)\) and \(y = (y_1, \ldots, y_d)\), the interlacing scheme orders the joint variable \(z = (z_1, \ldots, z_{2d})\) by \(z = (x_1, y_1, \ldots, x_{d}, y_d)\).
Typically, the variable \(x, y\) are ordered so that the correlation for \((x_{j}, x_{j'})\) and \((y_{j}, y_{j'})\) are strong only when the index \(j\) is close to \(j'\). Thus, capturing the model's spatial correlation structure means the ordering of the joint variable needs to preserve the proximity. Moreover, \(x_{j}\) and \(y_{j}\) represent the value of a high-dimensional stochastic process at the same site. Therefore, capturing the temporal correlation structure means the ordering of the joint variable needs to be such that \(x_j\) is placed close to \(y_j\). The proposed interlacing scheme preserves the spatial and temporal correlation structure as is desired.

\begin{figure}[t!]
    \centering
    \includegraphics[width = \textwidth]{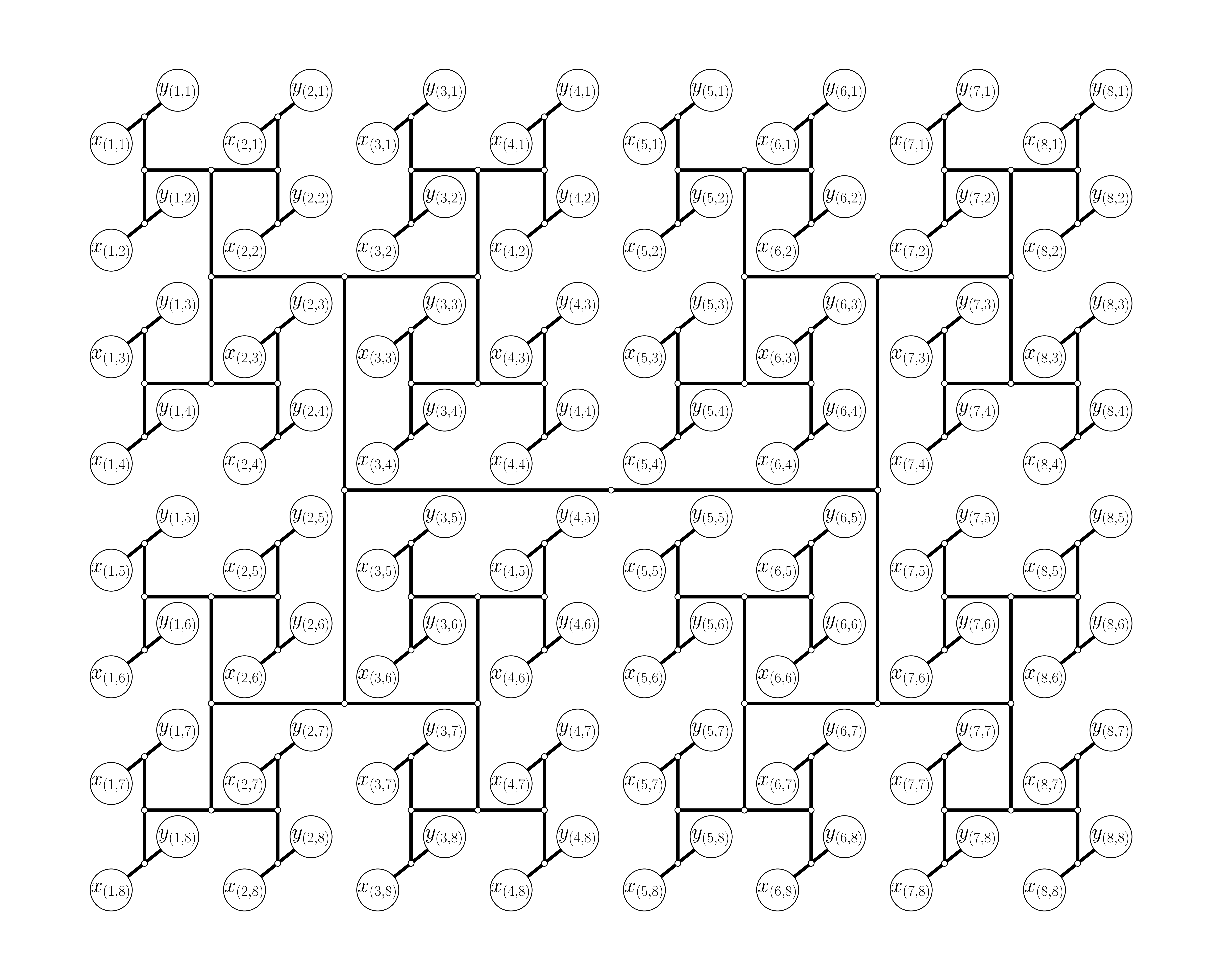}
    \caption{Illustration of the hierarchical bipartition of a 2D Ginzburg-Landau model discretized to an $8\times 8$ grid. The variable \(x = \left(x_{(i,j)}\right)_{1 \leq i, j\leq m}\) denotes the random variable at the initial time \(0\), and \(y = \left(y_{(i,j)}\right)_{1 \leq i, j\leq m}\) denotes the random variable at the intermediate time \(T-t\)}
    \label{fig:2D_GZ_bipar}
\end{figure}

An ordering of the variables leads to a corresponding hierarchical bipartition through \eqref{eqn: bipartition}. In Figure \ref{fig:2D_GZ_bipar}, we illustrate the resultant variable hierarchical bipartition in the 2D Ginzburg-Landau (G-L) model under the proposed interlacing scheme. The 2D G-L model has the natural geometry of a 2D grid, and we apply the natural construction of alternatively partitioning the \(x, y\) variables along the two axes of the 2D grid. One can see that the proposed bipartition structure for the joint variable indeed captures important correlation structures of the problem.

\subsection{Functional basis for the FTN ansatz}
The choice of the functional basis is intricately linked to our choice of the bounded domain \(x \in \mathcal{X}\). In this work, we choose a simple hypercube structure for the domain, i.e., \(\mathcal{X} = \prod_{j \in [d]}[a_j, b_j]\). This construction is beneficial in three important aspects. First, any bounded region one is interested in can always be covered by a hypercube, and so any domain \(\mathcal{X}\) can always be reduced to the hypercube case. Second, the simplicity in \(\mathcal{X}\) ensures the initial distribution \(p_0\) has a simple sampler during particle simulation. Most importantly, the separable structure of \(\mathcal{X}\) avoids the occurrence of the failure mode where spurious correlation from the irregularity of the domain prevents the joint density for \((X_0, X_{t})\) from having a low-rank structure. 

The choice of functional basis is a direct result of the hypercube structure \(\mathcal{X} = \prod_{j \in [d]}[a_j, b_j]\). As a minimal requirement, the chosen functional basis for the \(x\) variable needs to allow the functional hierarchical tensor to faithfully represent the uniform distribution on \(\mathcal{X}\). Thus, for each \(x_j\) variable, we pick the functional basis to be a collection of the orthonormal basis of \(L^{2}\left([a_j, b_j]\right)\) with the mild requirement that its span contains the constant function, i.e., the indicator function \(\chi_{[a_j, b_{j}]}\).

In practice, it is cumbersome to maintain a separate function basis for each variable. For our implementation, we use a fixed orthonormal basis \(\{\psi_{i}\}_{i = 1}^{n} \subset L^{2}\left([-1, 1]\right)\) with \(\psi_1 = \frac{1}{\sqrt{2}}\chi_{[-1, 1]}\). The function basis for \(x_j\) can be constructed by applying an appropriate linear change of co-ordinate to \(\{\psi_{i}\}_{i = 1}^{n}\) to obtain a basis in \(L^{2}\left([a_j, b_j]\right)\). In the Ginzburg-Landau model, the distribution of \(X_{t}\) is effectively compactly supported inside the hypercube \([-2, 2]^d\), which is why the functional basis for the \(y\) variable is also obtained by change of co-ordinate to \(\{\psi_{i}\}_{i = 1}^{n}\).

\subsection{Functional hierarchical sketching for joint density}
As the construction for the variable hierarchical bipartition structure and the basis are fixed, one can apply the functional hierarchical tensor sketching algorithm in \cite{tang2023solving} to obtain an FHT-based characterization of the joint density function for \((X_0, X_{t})\). The details of the sketching algorithm can be found in \cite{tang2023solving}.

We remark that the generalization to a variable-dependent function basis requires only a straightforward modification to the subroutine. 
Let \( \left\{ w^{(i)}:=\left(w_1^{(i)}, \ldots, w_{2d}^{(i)} \right) \right\}_{i = 1}^{N}\) be the samples of the joint variable \((X_0, X_{t})\) ordered according to the proposed interlacing scheme. In our case, we choose to obtain the variable-dependent functional basis \(\{\psi_{i;j}\}_{i = 1}^{n}\) through the aforementioned linear change of co-ordinate on a fixed basis function \(\{\psi_{i}\}_{i = 1}^{n} \subset L^{2}\left([-1, 1]\right)\). Thus, it is mathematically equivalent to apply a change of co-ordinate to the received sample \( \left\{ w^{(i)} \right\}_{i = 1}^{N}\) and then perform the sketch-based density estimation algorithm with the variable-invariant functional basis \(\{\psi_{i}\}_{i = 1}^{n}\), which is the setting in \cite{tang2023solving}. In general, to accommodate the variable-dependent function basis, the modified sketching procedure would use \(\{\psi_{i;j}\}_{i = 1}^{n}\) for the functional basis evaluation on the \(j\)-th co-ordinate of each sample. 

\begin{figure}[!ht]
    \centering
    \includegraphics[width = 0.6\textwidth]{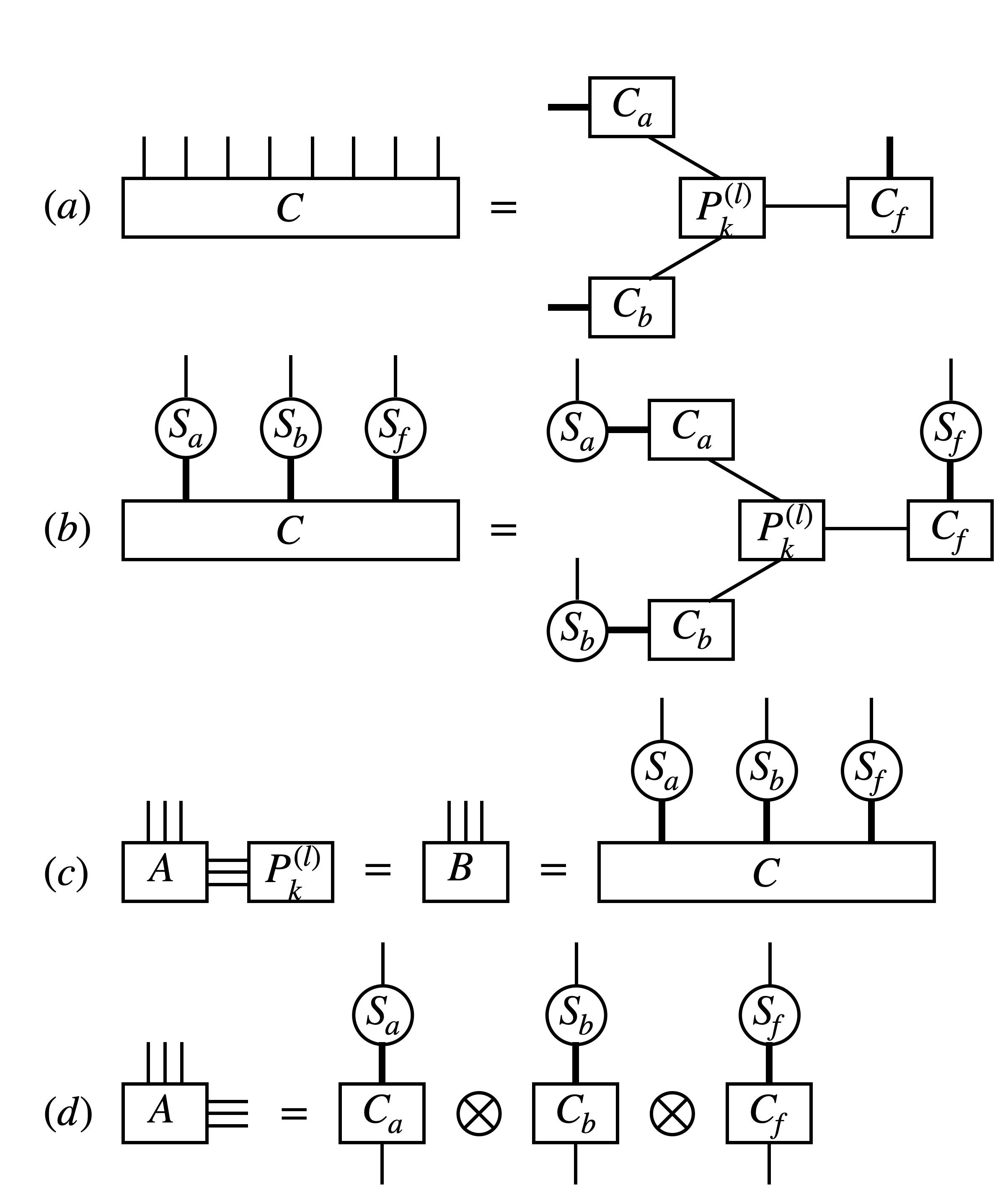}
    \caption{Tensor diagram representation of main equations used in the hierarchical sketching algorithm. Equation \eqref{eqn: unsketched linear system for G} is summarized in (a). Equation \eqref{eqn: sketched linear system for G} is shown in (b) and (c). The coefficient term in \eqref{eqn: sketched linear system for G} is shown in (d).}
    \label{fig:tensor_diagram_equation}
\end{figure}

To give the reader some intuition on the approach, we give a condensed description of the main procedure underlying the sketch-based algorithm for the \(2d\)-dimensional joint density function for \((X_0, X_{t})\). For the reader's convenience, relevant equations are included in Figure \ref{fig:tensor_diagram_equation} in terms of the tensor diagram.

Let the tensor cores for the joint density \(g\) be denoted by \(\{P_{k}^{(l)}\}_{k \in [2^l], l = 0, \ldots, L+1}\), where \(P_{k}^{(l)}\) is the \(k\)-th tensor core at level \(l\).
For simplicity, we outline the procedure for solving for the tensor core at intermediate levels with $0< l < L+1$. 

Using the notation for hierarchical bipartition as in \eqref{eqn: bipartition}, write \(a := I_{2k-1}^{(l+1)}, b := I_{2k}^{(l+1)}\) and \(f = [2d] - a \cup b\).
In this case, \(P_{k}^{(l)}\) is a \(3\)-tensor \(P_{k}^{(l)}\colon [r_a] \times [r_b] \times [r_f] \to \R\) for some \(r_a, r_b, r_f \in \mathbb{N}\).
The structural low-rankness property of \(g\) in \eqref{eqn: htn forward map} implies that there exist \(C_a \colon [n^{|a|}] \times [r_a] \to \R\), \(C_b \colon [n^{|b|}] \times [r_b] \to \R\), \(C_f \colon [r_f] \times [n^{|f|}] \to \R\) such that the following equation holds (illustrated in Figure \ref{fig:tensor_diagram_equation}(a)):
\begin{equation}\label{eqn: unsketched linear system for G}
    C(i_{1}, \ldots, i_{2d}) = \sum_{\alpha, \beta, \theta}C_{a}(i_{a}, \alpha)C_{b}(i_{b}, \beta)P_{k}^{(l)}(\alpha, \beta, \theta)C_{f}(\theta, i_{f}).
\end{equation}

Equation \eqref{eqn: unsketched linear system for G} is exponential-sized, and thus, one needs to use sketching to solve a system of smaller size.
The hierarchical sketching algorithm essentially solves the over-determined linear system \eqref{eqn: unsketched linear system for G} for \(P_{k}^{(l)}\) with the use of sketch functions. Let \(\tilde{r}_a, \tilde{r}_b, \tilde{r}_f\) be integers so that \(\tilde{r}_a > r_a, \tilde{r}_b > r_b, \tilde{r}_f > r_f\). Through the sketch functions \(S_a \colon [n^{|a|}] \times [\tilde{r}_a] \to \R, S_b \colon [n^{|b|}] \times [\tilde{r}_b] \to \R, S_f \colon [n^{|f|}] \times [\tilde{r}_f] \to \R\), one can contract the linear system in \eqref{eqn: unsketched linear system for G} with \(S_a, S_b, S_f\) at the variables \(i_a, i_b, i_f\) respectively, which leads to the following sketched linear system (illustrated in Figure \ref{fig:tensor_diagram_equation}(b)-(d)):
\begin{equation}\label{eqn: sketched linear system for G}
    B(\mu, \nu, \zeta) = \sum_{\alpha, \beta, \theta}A_{a}(\mu, \alpha)A_b(\nu, \beta)A_f(\zeta, \theta)P_{k}^{(l)}(\alpha, \beta, \theta),
\end{equation}
where \(A_{a}, A_b, A_f\) are respectively the contraction of \(C_a, C_b, C_f\) by \(S_a, S_b, S_f\) and \(B\) is the contraction of \(C\) by \(S_a \otimes S_b \otimes S_f\). The equation \eqref{eqn: sketched linear system for G} can be seen as the following linear system of \(\tilde{r}_a\tilde{r}_b\tilde{r}_f\) equations for the unknown \(P_{k}^{(l)}\):
\[
(A_{a} \otimes A_{b} \otimes A_{f})P_{k}^{(l)} = B.
\]

In practice, one has access to \(g\) through the received sample \(\left\{ w^{(i)} \right\}_{i = 1}^{N}\). Through somewhat involved derivations, which are outlined in \cite{tang2023solving}, one shows that \(B, A_{a}, A_b, A_f\) can be obtained by calculating statistical moments of the distribution \(g\), which can be done efficiently through sample-based moment estimation. Thus, the hierarchical sketching algorithm systematically collects statistical moments and subsequently solves tensor cores by solving independent linear systems.

\section{Solving for the Kolmogorov backward equations}\label{sec: alg}
In this section, we show how the Markov operator gives rise to an approximated solution to the Kolmogorov backward equation. Section \ref{sec: main workflow} summarizes the main workflow used to solve the Kolmogorov backward equation. Section \ref{sec: special f} discusses solving for the Kolmogorov backward equation when the terminal condition admits an FHT ansatz and comments on the case for the general terminal condition.
Section \ref{sec: extension} discusses solving for the Kolmogorov forward equation and solving Kolmogorov equations under general SDE dynamics.

\subsection{Main workflow}\label{sec: main workflow}
We summarize our approach to solving the Kolmogorov backward equation. The given information regarding the PDE includes the potential function \(V \colon \mathbb{R}^{d} \to \mathbb{R}\), a hypercube region \(\mathcal{X}\), and terminal time \(T\). The main workflow is summarized in the following simple steps:
\begin{enumerate}
    \item Select K time steps \(0 = t_0 \leq t_1 \leq \ldots \leq t_{K} = T\). Then, with a chosen numerical scheme (e.,g. Euler-Maruyama), run \(N\) independent instance of SDE simulations on the stochastic dynamic equation \eqref{eqn: langevin} to time \(T\) with the initial condition \(X_{0} \sim \mathrm{Unif}(\mathcal{X})\). Obtain \(N\) trajectories and record the trajectory data \(\{X^{(i)}_{t_{j}}\}_{j \in [K], i \in [N]}\).

    \item For \(j = 1, \ldots, K\), obtain an approximated joint density \(g_j(x, y) \approx \mathbb{P}\left[X_{0} = x, X_{t_j} = y\right]\) by calling the sketch-based density estimation algorithm on the sample \(\left\{\left(X^{(i)}_{0}, X^{(i)}_{ t_j}\right)\right\}_{i \in [N]}\). Obtain the Markov operator \(G_{t}\) at \(t = t_j\) through \(G_{t_j}(x, y) = \mathrm{Vol}(\mathcal{X})g_j(x, y)\).
    
    \item For any \((x, t)\), output the continuous-in-time solution for \(u(x, t)\), defined as follows: Suppose that \(t \in [0, T]\) satisfies \(T-t \in (t_j, t_{j+1})\), then \(u(x, t)\) can be approximated with \(u_{j}(x)\) and \(u_{j+1}(x)\) via an appropriate interpolation. The values of \(u_{j}(x)\) and \(u_{j+1}(x)\) are obtained with the corresponding Markov operator. The details are in Section \ref{sec: special f}.
\end{enumerate}

\subsection{FHT ansatz for the solution}\label{sec: special f}

After approximating the Markov operator \(G_{t}(x, y)\) on \(x \in \mathcal{X}\) via the sketch-based density estimation, one obtains the value of \(u(x, T-t)\) for \(x \in \mathcal{X}\) through the following equation. Let \(g_{\mathrm{FHT}}\) denote the functional hierarchical tensor ansatz for the joint density obtained by the sketch-based density estimation algorithm. One then has
\begin{equation}\label{eqn: FHT form of u}
    u(x, T-t) = \int_{\Omega} G_{t}(x, y)f(y) \, dy  \approx \mathrm{Vol}(\mathcal{X})\int_{\Omega} g_{\mathrm{FHT}}(x, y)f(y) \, dy, \quad x \in \mathcal{X}.
\end{equation}
In this subsection, we discuss the special case in which the terminal condition \(f\) admits an FHT ansatz. In this case, we can indeed show that the solution \(u(x, T-t)\) admits an FHT ansatz through an efficient tensor contraction. 

We give two motivating examples in which \(f\) admits an FHT ansatz. First, when \(f\) is separable, i.e. when \(f(y) = \prod_{j \in [d]}h_{j}(y_j)\), one can show that \(f\) is of rank one across any hierarchical bipartition. Second, when \(f\) is a sum of univariate functions, i.e., when \(f(y) = \sum_{j \in [d]}h_{j}(y_j)\), one can show that \(f\) is of rank two across any hierarchical bipartition. For general terminal condition \(f\), one can likewise apply an FHT-based approximation to the terminal condition. Discussion on the approximation of general functions by functional hierarchical tensors can be found in \cite{bachmayr2016tensor, ballani2013black} and references therein. At the end of this subsection, we give the procedure for obtaining the solution for general \(f\) without FHT approximations.

We now construct the FHT ansatz for the PDE solution \(u \approx \mathrm{Vol}(\mathcal{X})\int_{\Omega} g_{\mathrm{FHT}}(x, y)f(y) \, dy\). Suppose \(d = 2^{L}\) where \(L \geq 2\). We use the notational convention that \(G_{k}^{(l)}\) denotes the tensor core of an FHT ansatz corresponding to the \(k\)-th node at the \(l\)-th level. When the node is not a root nor leaf node, the tensor core \(G_{k}^{(l)}\) is a \(3\)-tensor and is otherwise a \(2\)-tensor. The tensor indexing convention strictly follows that of \cite{tang2023solving}.

\begin{figure}[!h]
    \centering
    \includegraphics[width = 0.9\textwidth]{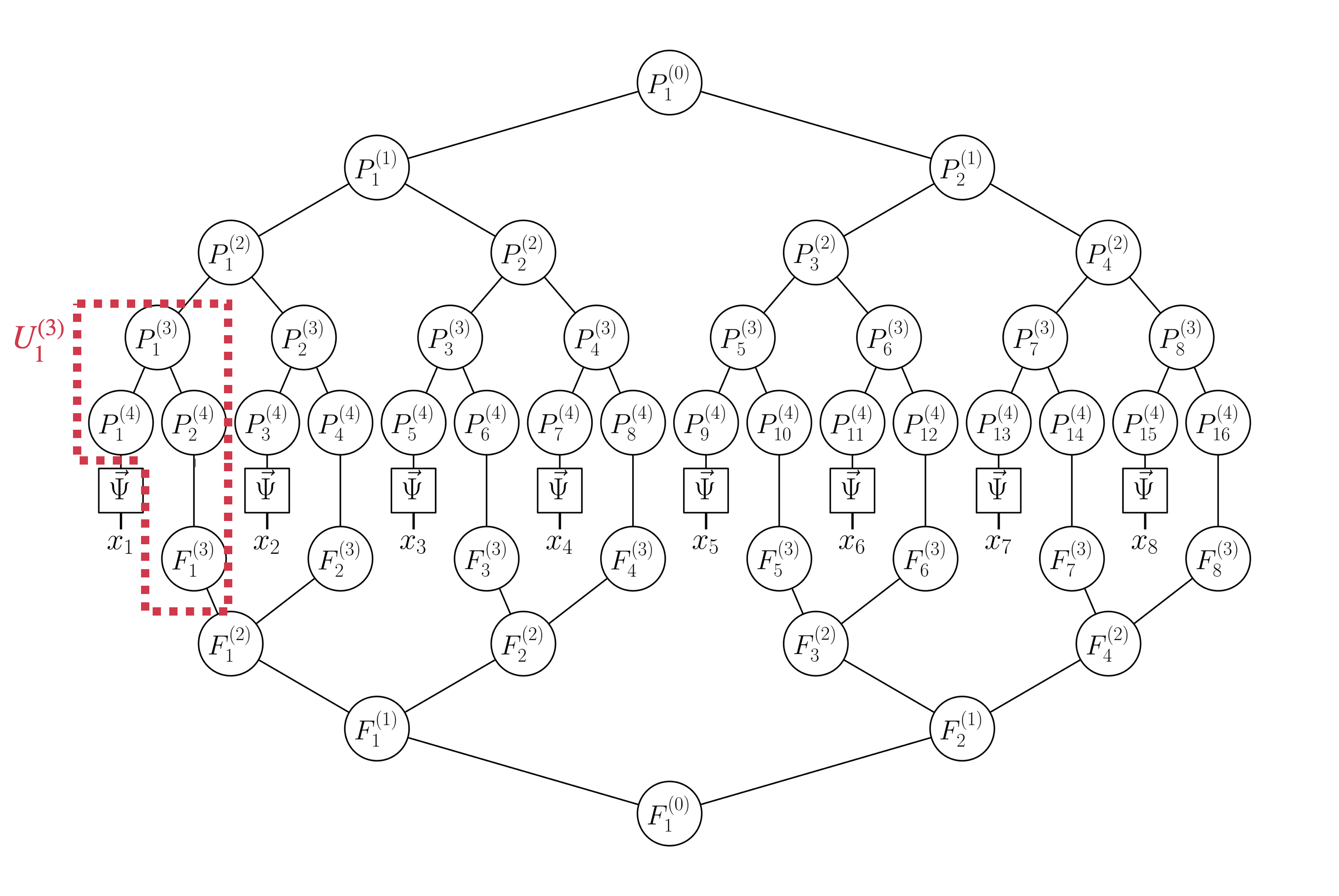}
    \caption{Illustration of the tensor contraction diagram for the approximated solution to the Kolmogorov backward equation for \(d = 8\).}
    \label{fig:FHT_contraction}
    \includegraphics[width = 0.9\textwidth]{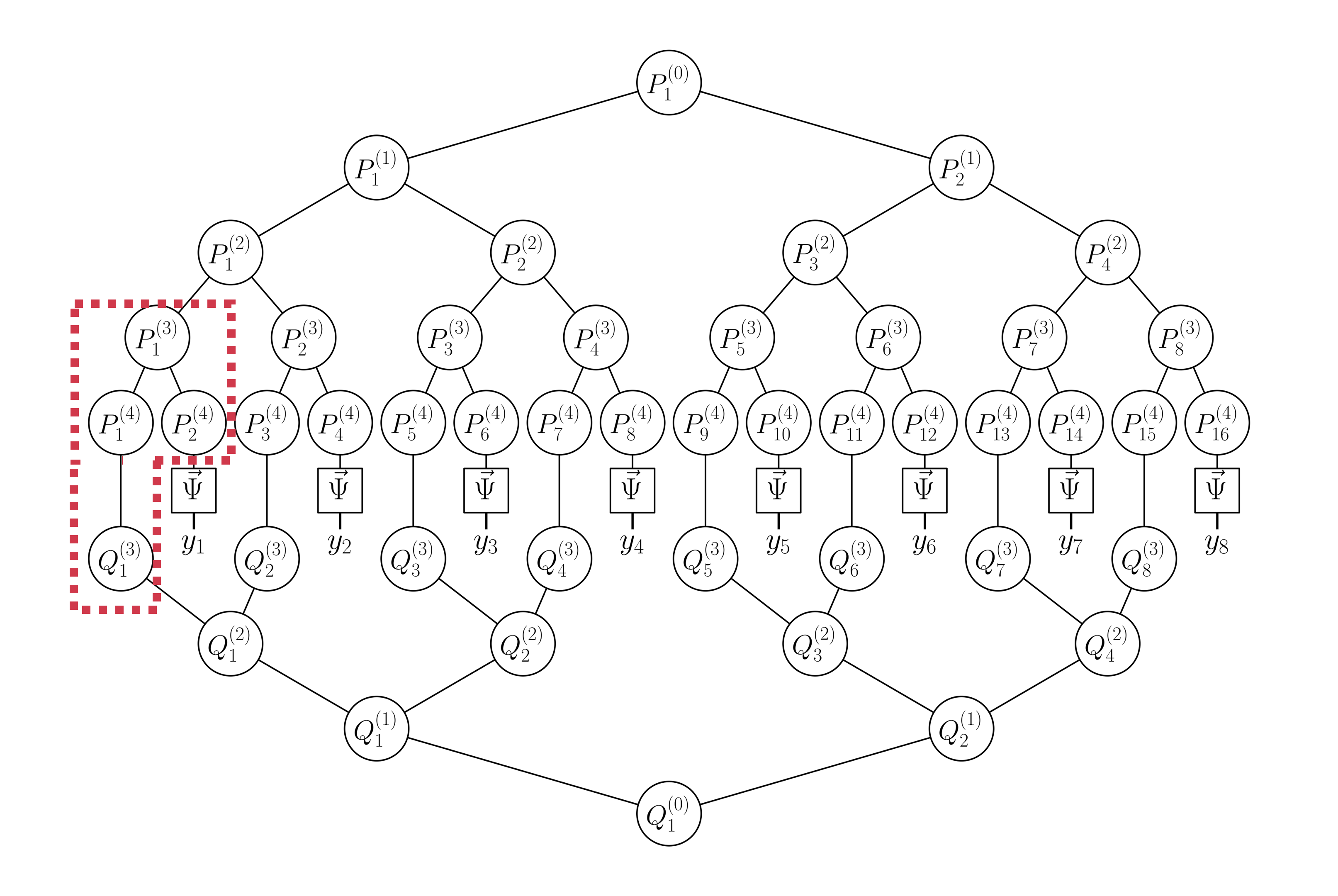}
    \caption{Illustration of the tensor contraction diagram for the approximated solution to the Kolmogorov forward equation for \(d = 8\).}
    \label{fig:FHT_contraction_forward}
\end{figure}

Let the tensor core of the solved joint density \(g_{\mathrm{FHT}}\) be denoted by \(\{P_{k}^{(l)}\}_{k \in [2^l], l = 0, \ldots, L+1}\), and let the tensor core of the terminal condition \(f\) be denoted by \(\{F_{k}^{(l)}\}_{k \in [2^l], l = 0, \ldots, L}\). We shall give an explicit FHT representation of \(u\) by a collection of tensor cores \(\{U_{k}^{(l)}\}_{k \in [2^l], l = 0, \ldots, L}\). Indeed, as \(f\) and \(g_{\mathrm{FHT}}\) are of compatible network structure, obtaining the FHT representation for \(\int_{\Omega} g_{\mathrm{FHT}}(x, y)f(y) \, dy\) is done through a simple tensor contraction. See Figure \ref{fig:FHT_contraction} for a tensor diagram illustration of the procedure. 

As the illustrated tensor diagram indicates, to form the FHT ansatz for \(u\), the construction is through combining the tensor cores from \(g_{\mathrm{FHT}}\) and \(f\). For the root level \(l = 0\), we construct the tensor core \(U_{1}^{(0)}\) as the Kronecker product between two \(2\)-tensors \(P_{1}^{(0)}\) and \(F_{1}^{(0)}\), where in addition a coefficient term is applied to absorb the normalization constant from the volume of \(\mathcal{X}\):
\[
U_{1}^{(0)}\left((\alpha, \mu), (\beta, \nu)\right) = \mathrm{Vol}(\mathcal{X})P_{1}^{(0)}(\alpha, \beta)F_{1}^{(0)}(\mu, \nu).
\]

For an intermediate level \(l \not \in \{ 0, L\}\) and an arbitrary block index \(k \in [2^l]\), the tensor core \(U_{k}^{(l)}\) is the \(3\)-tensor formed by a tensor-equivalent of Kronecker product between the two \(3\)-tensors \(P_{k}^{(l)}\) and \(F_{k}^{(l)}\):
\[
U_{k}^{(l)}\left((\alpha, \mu), (\beta, \nu), (\theta, \xi)\right) = P_{k}^{(l)}(\alpha, \beta, \theta)F_{k}^{(l)}(\mu, \nu, \xi).
\]

The construction for the leaf level is more involved, as one needs to perform tensor contraction over four tensor cores. In this case, for \(k \in [d]\), one constructs the \(2\)-tensor \(U_{k}^{(L)}\) by the following formula (see illustration in Figure \ref{fig:FHT_contraction}):
\[
U_{k}^{(L)}\left(i, (\theta, \xi)\right) = \sum_{i', \alpha, \beta} 
P_{2k-1}^{(L+1)}(i, \alpha)P_{2k}^{(L+1)}(i', \beta)P_{k}^{(L)}(\alpha, \beta, \theta)F_{k}^{(L)}(i', \xi).
\]

We remark that the approach likewise generalizes to the situation when terminal condition \(f\) is a sum of functions admitting FHT ansatz. In that case, one can contract each FHT with the FHT approximation of the Markov operator, which outputs the approximated PDE solution as a sum of functions with the FHT ansatz.

For the case of a general terminal condition, one can also obtain an approximated PDE solution. In this case, one can use the fact that conditional distribution with a functional hierarchical tensor is efficient with only an \(O(d)\) scaling. Thus, for \(x \in \mathcal{X}\), one uses \(g_{\mathrm{FHT}}\) to perform an approximate conditional sampling of \(X_{t}\) condition on \(X_0 = x\). After \(M\) samples \((y^{(i)})_{i =1}^{M}\) are collected, the solution to the PDE at \(u(x, T-t)\) can be approximated by 
\[
u(x, T-t) \approx \frac{1}{M}\sum_{i = 1}^{M}f(y^{(i)}).
\]

In this case, the FHT-based approximated solution is more efficient than simulating the conditional distribution of \(X_{t}\) based on SDE simulation over the dynamics in \eqref{eqn: langevin}. The reason is that the cost of calculating the drift term in the SDE simulation is proportional to the cost of evaluating the gradient of \(V\), as well as the total number of times in which the gradient is called. In practical cases, even in the case where the cost of evaluating \(\nabla V\) is \(O(d)\) and \(t = O(1)\), one would still keep a small step size for stability, which leads to a cost often much higher than the conditional sampling subroutine for FHT. Thus, the proposed approach can be quite beneficial when the cost of SDE simulation is much higher than the cost of conditional sampling.

\subsection{Extensions}\label{sec: extension}

\paragraph{FHT ansatz for solution to the Kolmogorov forward equation}

After approximating the Markov operator \(G_{t}(x, y)\) on \(x \in \mathcal{X}\) via the sketch-based density estimation, one likewise obtains the solution to the Kolmogorov forward equation through tensor contraction. Assume that \(\mathcal{X}\) is sufficiently large such that the support of the initial distribution \(q(x)\) is concentrated in \(\mathcal{X}\) except for a negligible volume. Then, the functional hierarchical tensor ansatz for the joint density outputs
\[
p(y, t) = \int_{\Omega} G_{t}(x, y)q(x) \, dx  \approx \mathrm{Vol}(\mathcal{X})\int_{\Omega} g_{\mathrm{FHT}}(x, y)q(x) \, dx, \quad x \in \mathcal{X},
\]
which is an approximate solution to the Kolmogorov forward equation with \(q\) as the initial distribution. We discuss the case in which \(q\) is in the form of an FHT ansatz or when \(q\) is characterized in terms of samples following the distribution governed by \(q\). The latter case can be converted to the former case by using the sketch-based algorithm to obtain an FHT approximation to \(q\) via the samples. Let the tensor core of \(g_{\mathrm{FHT}}\) be denoted by \(\{P_{k}^{(l)}\}_{k \in [2^l], l = 0, \ldots, L+1}\), and let the tensor core of \(q\) be denoted by \(\{Q_{k}^{(l)}\}_{k \in [2^l], l = 0, \ldots, L}\). The procedure of deriving the tensor core for the FHT approximation of \(p\) follows exactly as the tensor contraction in Section \ref{sec: special f}, and thus we omit the details. The procedure is illustrated via a tensor diagram in Figure \ref{fig:FHT_contraction_forward}

\paragraph{Kolmogorov equations beyond Langevin dynamics}
Though we have been focusing on Langevin dynamics for simplicity, our method can be generalized to a wider range of SDE dynamics. In particular, one can consider a general particle system with the following form:
\begin{equation}\label{eqn: general sde for fokker-planck}
    dX_t = f(X_t, t) dt + \sigma(X_t, t)dB_{t},
\end{equation}
where \(f \colon \R^d \times [0, T] \to \R^d\) and \(\sigma \colon \R^d \times [0, T] \to \R^{d\times d}\) are general advection and diffusion terms. The formula for the Kolmogorov forward equation is
\begin{equation*}
    \partial_t p = \nabla \cdot (p f 
) + \frac{1}{2}\sum_{i, j, k =1}^{d} \partial_{x_i x_{j}} (\sigma_{ik}\sigma_{jk} p), \quad x \in \Omega \subset \R^{d}, t \in [0, T].
\end{equation*}
In this case, solving for the forward equation requires only simple modifications to the proposed routine. In this case, for a hypercube region \(\mathcal{X}\), one can sample \(X_0 \sim \mathrm{Unif}(\mathcal{X})\) and then perform the SDE simulation with \eqref{eqn: general sde for fokker-planck}. The Markov operator \(G_{t}\) can likewise be obtained through the probability density function via 
\begin{equation*}
    G_{t}(x, y) = \mathrm{Vol}(\mathcal{X})\mathbb{P}_{X_0 \sim \mathrm{Unif}(\mathcal{X})}\left[X_{0} = x, X_{t} = y \right], \quad x \in \mathcal{X},
\end{equation*}
and the solution \(p(y, t)\) can be obtained by contracting the Markov operator with the prescribed initial condition \(q\).

The situation for the Kolmogorov backward equation is more involved. In this case, the PDE is
\[
\partial_t u + \nabla u \cdot  f 
 + \frac{1}{2}\sum_{i, j, k =1}^{d}  \sigma_{ik}\sigma_{jk} \partial_{x_i x_{j}}u = 0, \quad x \in \Omega \subset \R^{d}, t \in [0, T].
\]

Importantly, the dynamics \eqref{eqn: general sde for fokker-planck} is no longer time-invariant, and the relevant solution operator admits the representation
\[
G_{t,T}(x, y) := \mathbb{P}\left[X_{T} = y \mid X_t = x \right] = \frac{1}{p_0(x)}\mathbb{P}\left[X_{t} = x, X_{T} = y \right].
\]
In this case, obtaining the value of \(u(x, T-t)\) at each \(t \in [0, T]\) requires a separate set of SDE simulations. For a hypercube region \(\mathcal{X}\), one can sample \(X_t \sim \mathrm{Unif}(\mathcal{X})\) and then perform the SDE simulation with \eqref{eqn: general sde for fokker-planck}. The solution operator \(G_{t,T}\) can likewise be obtained through the probability density function via 
\begin{equation*}
    G_{t,T}(x, y) = \mathrm{Vol}(\mathcal{X})\mathbb{P}_{X_t \sim \mathrm{Unif}(\mathcal{X})}\left[X_{t} = x, X_{T} = y \right], \quad x \in \mathcal{X}.
\end{equation*}
Once the solution operator is obtained, the solution can be similarly obtained either through conditional sampling or tensor contraction, depending on whether \(f\) admits an FHT-based representation.









\section{Application to the Ginzburg-Landau model}\label{sec: numerics}
We apply the proposed Markov operator method to a discretized Ginzburg-Landau model in 1D and 2D. For the sake of simplicity, we solve for the backward equation at \(t = 0\). In general, solving for the backward equation solution at \(t \in [0, T]\) can be done by running the algorithm at multiple time slots as prescribed in Section \ref{sec: alg}.

\subsection{Implementation detail}

\paragraph{The choice of basis function}
For the Ginzburg-Landau model, the $\int_{a} |1 - x(a)^2|^2 da$ term ensures that the distribution of each entry of \(X_{T}\) is essentially bounded within the interval \([-2, 2]\). For simplicity, we consider the situation whereby one is interested in the value of the KBE solution \(u\) for \(\mathcal{X} = [-2, 2]^d\). In the proposed method for the Markov operator estimation procedure, the sample is obtained by generating \(X_0 \sim \mathrm{Unif}([-2, 2]^d)\) and subsequently generating \(X_{t}\) by simulating the SDE in \eqref{eqn: langevin}. Thus, the joint density function for the variable \((X_0, X_{t})\) can be modelled to have compact support in \([-2, 2]^{2d}\). Therefore, we use the polynomial representation choosing a maximal degree parameter \(q\) and picking the first \(n = q + 1\) Legendre basis polynomial in $[-2,2]$, denoted \(\{\psi_{i}\}_{i=0}^{q}\).

\paragraph{Bipartition structure}
Similar to the hierarchical bipartition proposed in \cite{tang2023solving}, we recursively bipartition the Cartesian grid along its axes in an alternating fashion. The structure for the two-dimensional case is shown in Figure \ref{fig:2D_GZ_bipar} and can be done likewise to other dimensions. Explicitly, assume for simplicity that a function \(u\) in $\Delta$-dimensions is discretized to the \(d=m^\Delta\) points denoted \(\{u_{i_1, \ldots, i_\Delta}\}_{i_1, \ldots, i_\Delta \in [m]}\) with \(m = 2^{\mu}\). In the Markov operator case, one considers the discretization of \(u\) at initial and terminal time. Thus, one defines variables \(\{x_{i_1, \ldots, i_\Delta}, y_{i_1, \ldots, i_\Delta}\}_{i_1, \ldots, i_\Delta \in [m]}\), where \(x_{i_1, \ldots, i_\Delta}, y_{i_1, \ldots, i_\Delta}\) respectively corresponds to the discretization of \(u\) at the initial and terminal time for the location \((i_1, \ldots, i_\Delta)\). For each index \(i_\delta\), one performs a length \(\mu\) binary expansion $i_\delta = a_{\delta 1}\ldots a_{\delta \mu}$.
The variable \(x_{i_1, \ldots, i_\Delta}, y_{i_1, \ldots, i_\Delta}\) is respectively given index \(2k - 1, 2k\) where \(k\) has the length \(\mu \Delta\) binary expansion \(k = a_{11} a_{21} \ldots a_{\Delta-1,\mu} a_{\Delta,\mu}\). The associated binary decomposition structure follows from the given index according to \eqref{eqn: bipartition}. One can check that the construction exactly corresponds to performing binary partition by sweeping along the $\Delta$ dimensions. For intuition on this construction, one can check that the bipartition in Figure \ref{fig:2D_GZ_bipar} corresponds to the indexing \(k = a_{11} a_{21} a_{12} a_{22} a_{13} a_{23}\).

\paragraph{The choice of sketch function}
For the density estimation procedure under hierarchical sketching, the quality of the sketch function is critical to the fidelity of the Markov operator approximation. 
In our work, the sketch functions we choose are of two kinds. The first kind of sketch function consists of monomials of the Legendre basis of low degree. For \(a = I_{k}^{(l)}\) or \(a = [2d] - I_{k}^{(l)}\) (see definition in \eqref{eqn: bipartition}), the chosen sketch functions are of the following form:
\[
f(z_a) = \prod_{j \in a}{\psi_{i_j}(z_j)},
\]
for some choice of $0 \le i_j \le q$ for each $j$. One typically chooses \(f\) with small total degree defined by \(\mathrm{deg}(f) = \sum_{j\in a} i_j\). 

The second kind of sketch function is motivated by the renormalization group, in particular by the coarse-graining of the variables. In particular, one chooses a cluster \(h \subset a\) and the resultant function \(g_{h}\) is the average over variables in \(h\):
\[
g_{h,i}(z_a) \equiv \frac{1}{|h|}\sum_{j \in h}z_j.
\]
One can augment the sketch function set by including terms such as \(g_{h,1}^2\). Empirically, including such coarse-grained functions greatly increases model performance and stability. 

\paragraph{Sketch postprocessing}
For efficient performance of the sketching operation, it is often advantageous to apply postprocessing to reduce the noise level for the sketching operation during its moment estimation tasks. In particular, an important subroutine of the tensor sketching operation boils down to performing the estimation of \(\mathbb{E}\left[g(W)\right]\) for a function \(g\) and a \(2d\)-dimensional random variable \(W = (X_0, X_{t})\).
Because the variables within \(X_0\) are pairwise independent, it turns out that one has \(\mathbb{E}\left[g(W)\right] = 0\) in a large number of cases. Therefore, we apply a simple hard truncation rule whereby we estimate \(\mathbb{E}\left[g(W)\right]\) to be zero when the sample mean for \(g(W)\) is not larger than twice the estimated standard deviation for the sample mean, and otherwise we take the sample mean as the estimate. While this step can be skipped when one uses a large number of samples for the FHT density estimation subroutine, we observe great practical noise reduction by applying this approach.

\paragraph{Solution postprocessing}
After the density estimation step is complete, one obtains the estimation of the joint density according to 
\[\mathbb{P}\left[X_0 = x, X_{t} = y\right] \approx g_{\mathrm{FHT}}(x, y),\]
and the Markov operator is estimated by
\(
G_t(x, y) \approx \mathrm{Vol}(\mathcal{X})g_{\mathrm{FHT}}(x, y).
\)
In practice, it does not hold that \(\int g_{\mathrm{FHT}}(x, y) \, dy = 1/\mathrm{Vol}(\mathcal{X})\), and one achieves higher accuracy by adopting the following approximation of \(G_t\):
\begin{equation}
    G_t(x, y) \approx \frac{g_{\mathrm{FHT}}(x, y)}{\int_{s}g_{\mathrm{FHT}}(x, s)}.
\end{equation}

Correspondingly, the approximation of the backward solution is obtained by
\begin{equation}\label{eqn: non-FHT form of u}
    u(x, T-t) = \int_{\Omega} G_{t}(x, y)f(y) \, dy  \approx \frac{1}{\int_{s}g_{\mathrm{FHT}}(x, s)}\int_{\Omega} g_{\mathrm{FHT}}(x, y)f(y) \, dy, \quad x \in \mathcal{X}.
\end{equation}

The formula in \eqref{eqn: non-FHT form of u} is a division between two FHT-based functions, which is why \eqref{eqn: non-FHT form of u} is not an FHT-based format for \(u(x, T-t)\). To avoid this shortcoming, one can perform a simple FHT-based interpolation algorithm to obtain an FHT-based representation for \(u\). When one has a function with a black-box access to its evaluation, the hierarchical sketching algorithm can be used as an interpolation algorithm, and we omit the details for the sake of simplicity.

\subsection{1D Ginzburg-Landau potential}\label{sec: 1D test}

We first consider a 1D Ginzburg-Landau model. The potential energy is defined as
\begin{equation}\label{eqn: 1D GZ model}
V(x_1, \ldots ,x_m) := \frac{\lambda}{2} \sum_{i=1}^{m+1}\left(\frac{x_{i} - x_{i - 1}}{h}\right)^2 + \frac{1}{4\lambda} \sum_{i = 1}^{m} \left(1 - x_i^2\right)^2,
\end{equation}
where \( h = \frac{1}{m+1} \) and \(x_0=x_{m+1}=0\). In particular, we fix \( m = 128 \), \( \lambda = 0.005 \) and \( \beta = \frac{1}{20} \). The dimension is \(d = 128\).

\begin{figure}[!ht]
  \centering
  \includegraphics[width = 0.8\textwidth]{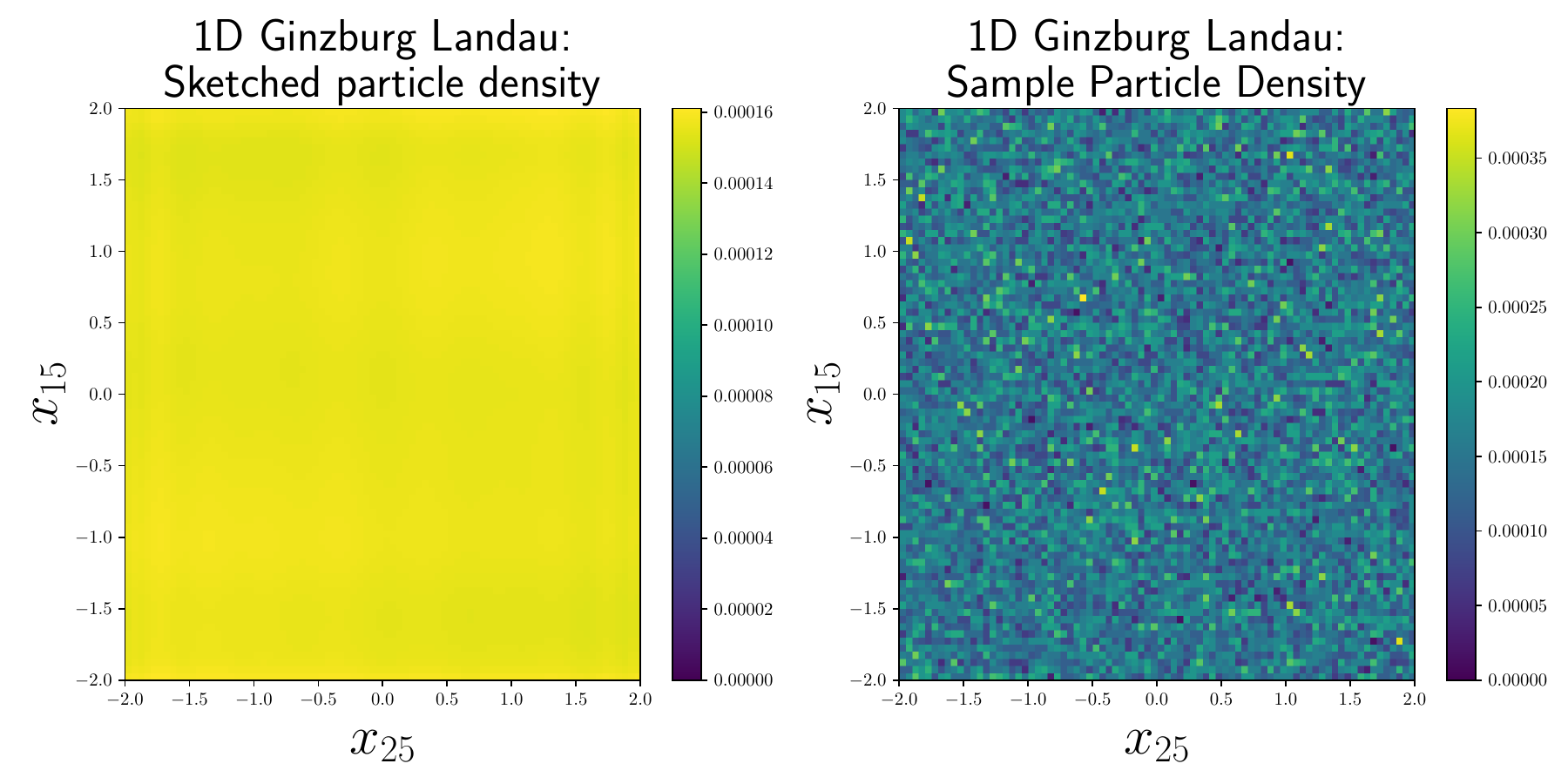}
\includegraphics[width = 0.8\textwidth]{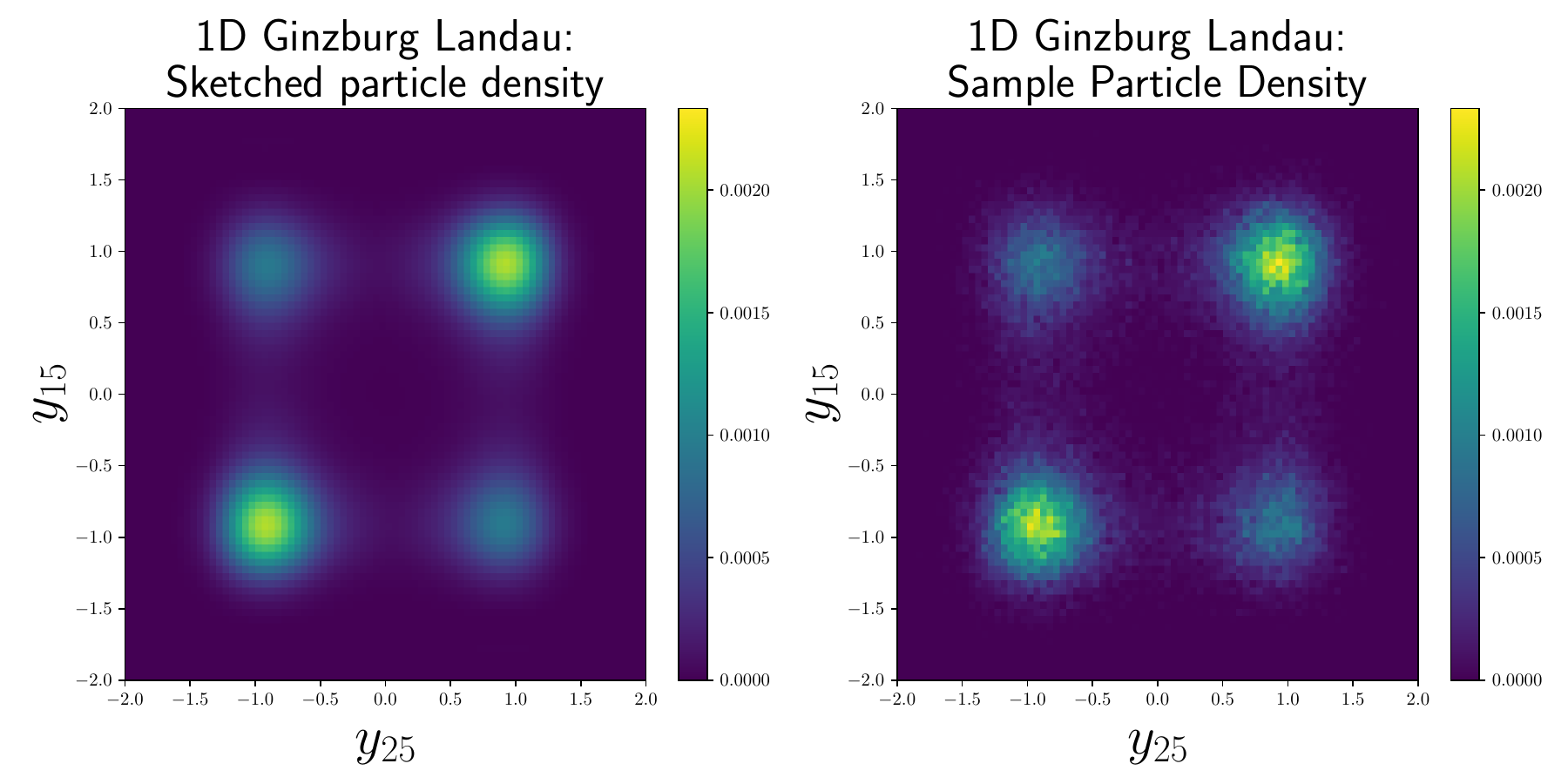}
  \caption{1D Ginzburg-Landau model. The plots of the marginal distribution at \((x_{15}, x_{25})\) and \((y_{15}, y_{25})\).}
  \label{Fig: 1D marginal}
\end{figure}

We perform \(N = 60000\) SDE simulations with \(T = 1\) and \( \delta t = \frac{1}{1000} \) with initial distribution $p_0 = \mathrm{Unif}([-2, 2]^{d})$. The maximal Legendre degree \(q\) is taken to be \(q = 20\), and the maximal internal bond dimension is \(r = 10\).

As an initial inquiry into the quality of the estimated joint density function \(\mathbb{P}\left[X_0 = x, X_T = y\right]\), Figure \ref{Fig: 1D marginal} compares the marginal distribution of \((x_{15}, x_{20})\) and \((y_{15}, y_{20})\) between empirical distributions and the model obtained by hierarchical tensor sketching. One can see that hierarchical tensor sketching can capture the independence between the \(x\) variables as well as the correlation between the \(y\) variables.

To assess the performance of the Markov operator in solving the backward equation, we propose to test the Markov operator with the following two terminal conditions:
\begin{equation}\label{eqn: terminal condition squared euclidean distance}
    f_{+}(y_1, \ldots, y_d) = \frac{1}{d}\sum_{j=1}^{d}(y_{j} - 1)^2, \quad f_{-}(y_1, \ldots, y_d) = \frac{1}{d}\sum_{j=1}^{d}(y_{j} + 1)^2,
\end{equation}
where the solution to the backward equation corresponds to the expected squared Euclidean distance from the two meta-stable states \(y_{+} = (1, \ldots, 1)\) and \(y_{-} = (-1, \ldots, -1)\) of the Ginzburg-Landau model.

\begin{figure}[!ht]
  \centering
  \includegraphics[width = 0.45\textwidth]{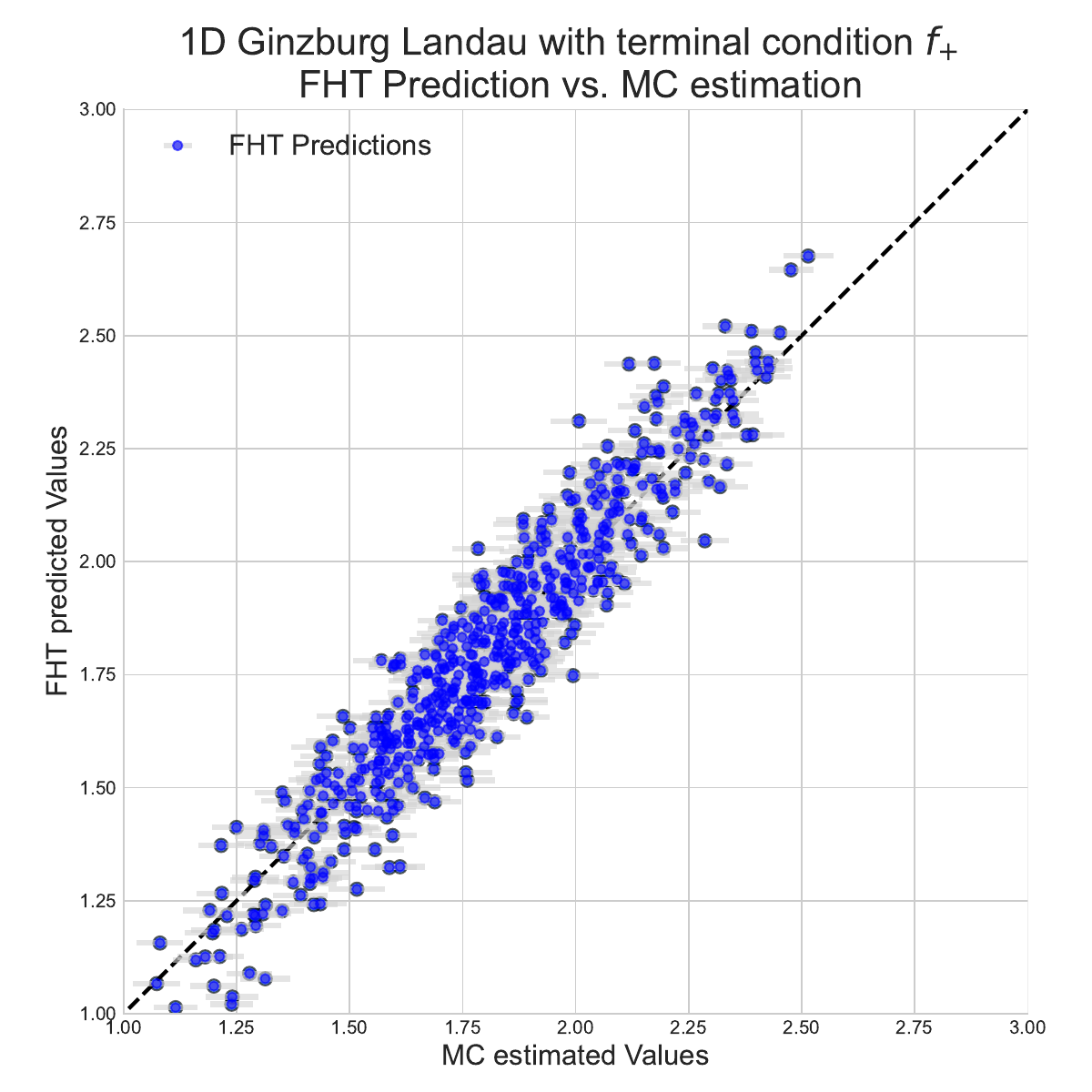}
  \includegraphics[width = 0.45\textwidth]{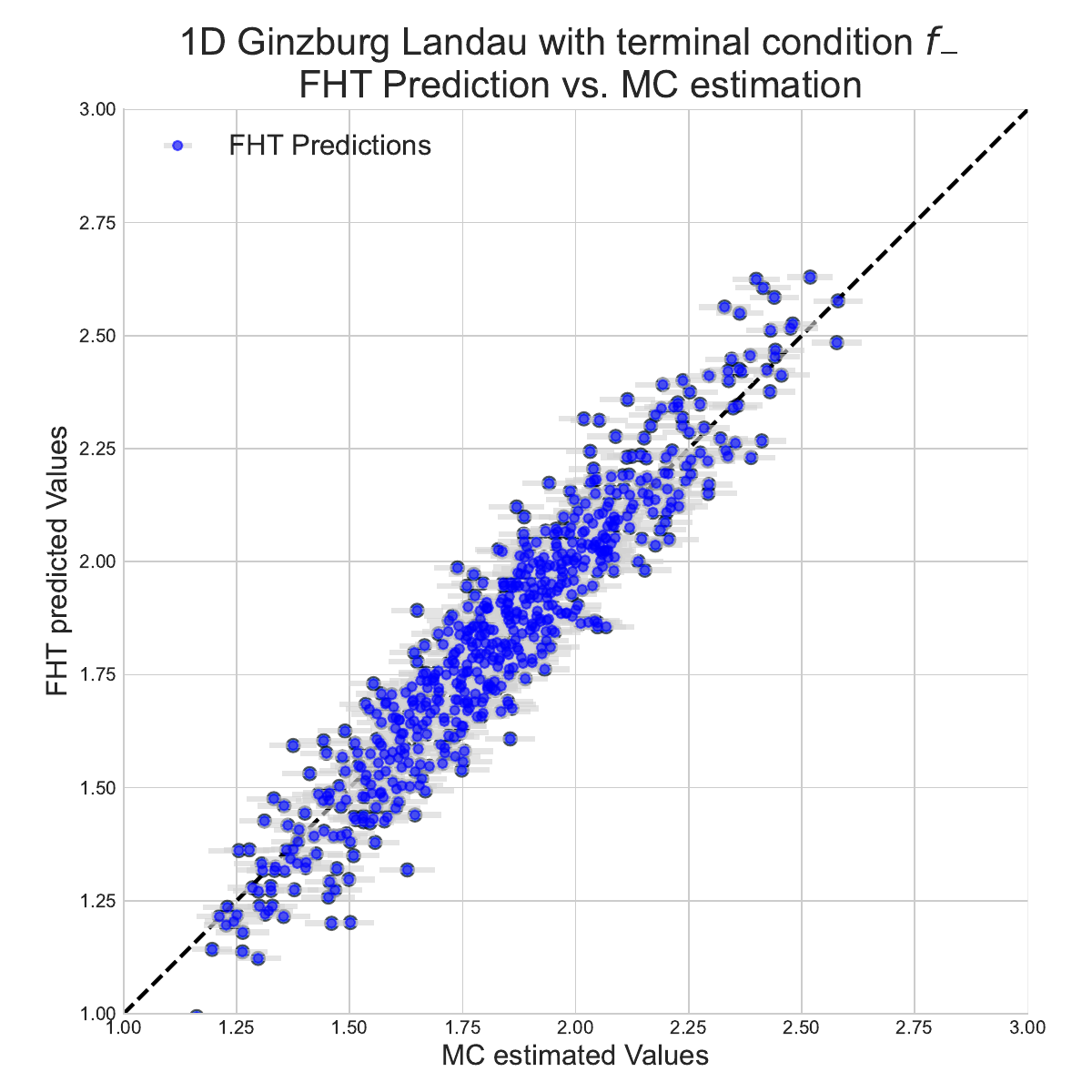}
  \caption{1D Ginzburg Landau model. The plot compares the function evaluation of the FHT-based approximation with reference values obtained from Monte-Carlo estimations, where the error bar indicates the standard deviation from each Monte-Carlo estimation. }
  \label{Fig: 1D error bar}
\end{figure}

For both terminal conditions, to compare the performance with the true backward equation, we estimate the value of \(u(x, 0)\) at \(N_{u} = 600\) randomly sampled points.
For each point, we sample \(100\) SDE trajectories to approximately query the evaluation of \(u\). The comparison with the predicted value through the Markov operator is plotted in Figure \ref{Fig: 1D error bar}, which shows that the approximated PDE solution matches the Monte Carlo estimation well. Overall, the average relative error for both terminal conditions is less than 4\%.

For the final numerical test, we consider using the Markov operator to probe the transitional properties of the Ginzburg-Landau model. In the Ginzburg-Landau model, one natural question to ask is whether a particular initial state \(x\) favors entering proximity to the state \(y_{+}\) or \(y_{-}\) under the Langevin dynamics. One way to analyze this is to consider two functions \(u_{+}, u_{-}\), which solves the backward equation for the following two terminal conditions:
\begin{equation}
    g_{+}(y_1, \ldots, y_d) = \exp{\left(-\frac{2}{d}\sum_{j=1}^{d}(y_{j} - 1)^2\right)}, \quad g_{-}(y_1, \ldots, y_d) = \exp{\left(-\frac{2}{d}\sum_{j=1}^{d}(y_{j} + 1)^2\right)},
\end{equation}
and the propensity for a state \(x\) to be near \(y_{+}\) is characterized the ratio \[\iota(x) = \frac{u_{+}(x, 0)}{u_{+}(x, 0) + u_{-}(x, 0)}.\] 

\begin{figure}[!ht]
  \centering
  \includegraphics[width = 0.8\textwidth]{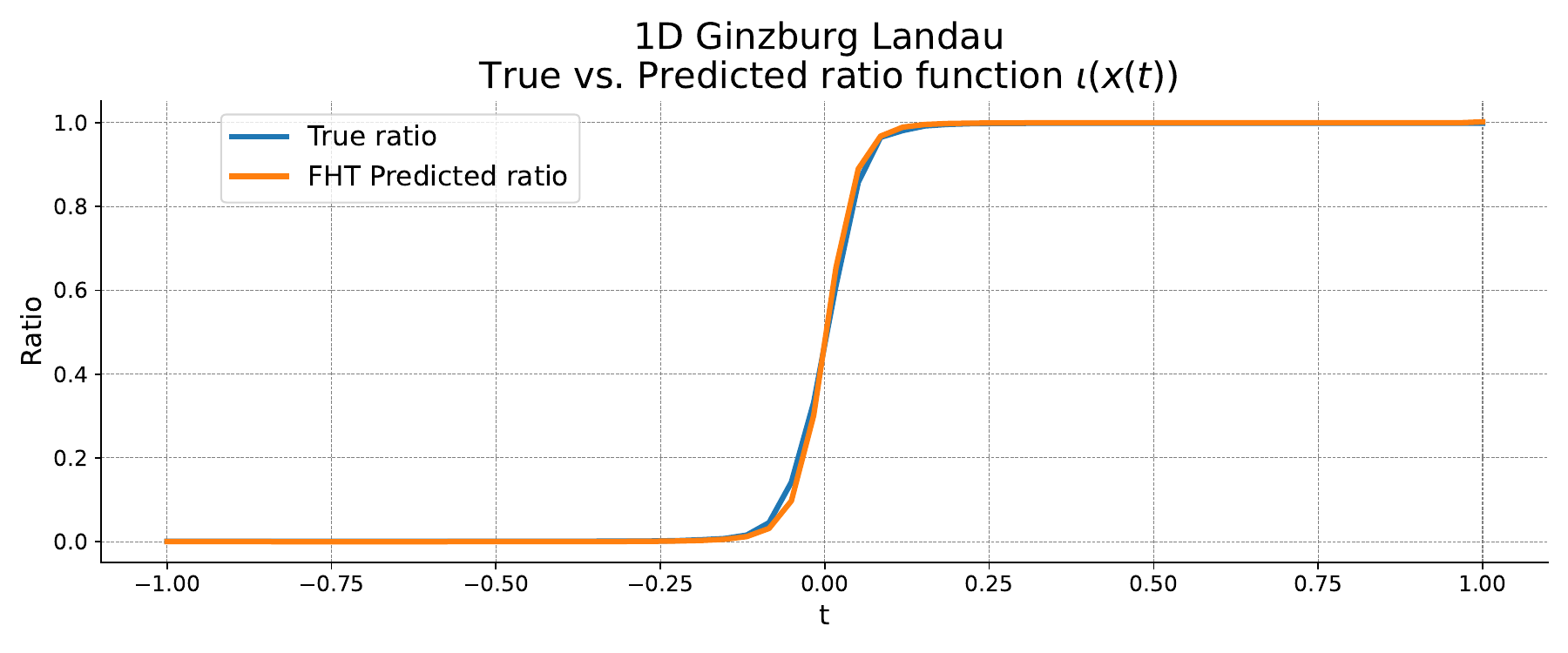}
  \caption{1D Ginzburg Landau model. Plot of \(\iota(x(t))\), i.e. the propensity of \(x(t) = (t,\ldots, t)\) to enter a proximity of \(y_{+} = (1, \ldots, 1)\). One can see that the Markov operator approach can capture the sharp transition at \(t = 0\).}
  \label{Fig: 1D ratio}
\end{figure}

In this case, both \(g_{+}\) and \(g_{-}\) are separable and admit a highly accurate FHT-based approximation, which is why one can directly solve for \(u_{+}, u_{-}\) through FHT contraction with the Markov operator. To consider an interesting case, we plot the predicted value of \(\iota(x)\) along the direction of \(x(t) = (t, \ldots, t) = t\mathbf{1}_{d}\). As can be seen in Figure \ref{Fig: 1D ratio}, \(\iota(x(t))\) has a sharp transition at the region of \(t \in [-0.15, 0.15]\), which indicates that the state \(x(t)\) has high likelihood to enter a proximity to \(y_{+}\) when \(t > 0\) is away from zero. Quite remarkably, this sharp transition is accurately captured by the FHT solution approximated by the Markov operator.

\subsection{2D Ginzburg-Landau potential}
We consider a 2D Ginzburg-Landau model. The potential energy is defined as
\begin{equation}
  V(x_{(1,1)}, \ldots, x_{(m,m)}) := \frac{\lambda}{2} \sum_{v \sim w}\left(\frac{x_{v} - x_{w}}{h}\right)^2 +  \frac{1}{4\lambda} \sum_{v}\left(1 - x_v^2\right)^2,
\end{equation}
where \(h = \frac{1}{m+1}\) and \( x_{(l, 0)} = x_{(l, m + 1)} = x_{(0, l)} = x_{(m+1, l)} = 0\) for \(l = 1, \ldots, m\). The parameters are \(m=16\), \(\lambda = 0.03\), \(\beta = \frac{1}{5}\). The dimension $d$ is $256$. 

We perform \(N = 60000\) SDE simulations with \(T = 1\) and \( \delta t = \frac{1}{1000} \), with the initial distribution $p_0 = \mathrm{Unif}([-2, 2]^{d})$. The maximal Legendre degree \(q\) is taken to be \(q = 10\), and the maximal internal bond dimension is \(r = 15\). 

\begin{figure}[!ht]
  \centering
  \includegraphics[width = 0.8\textwidth]{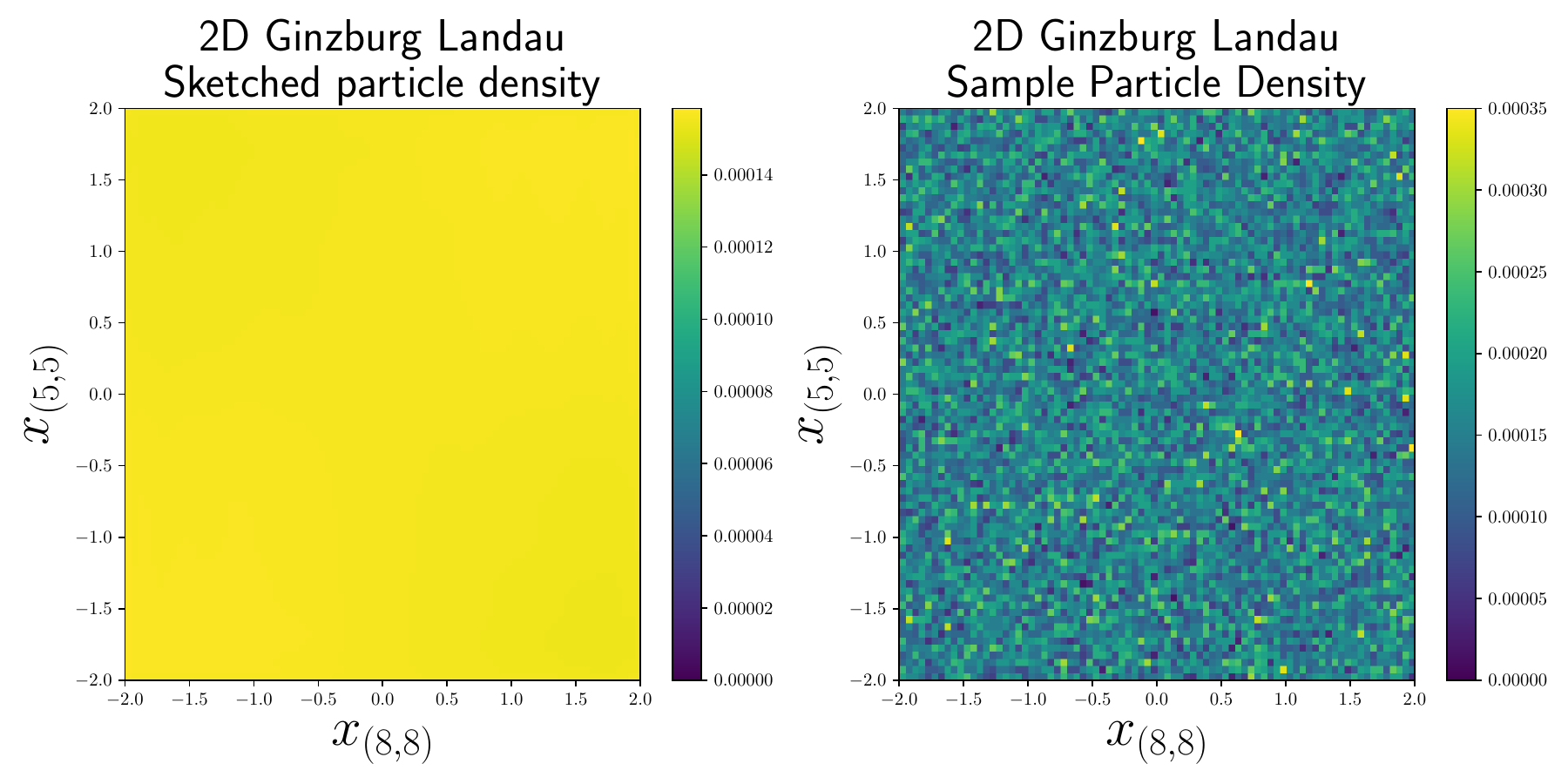}
\includegraphics[width = 0.8\textwidth]{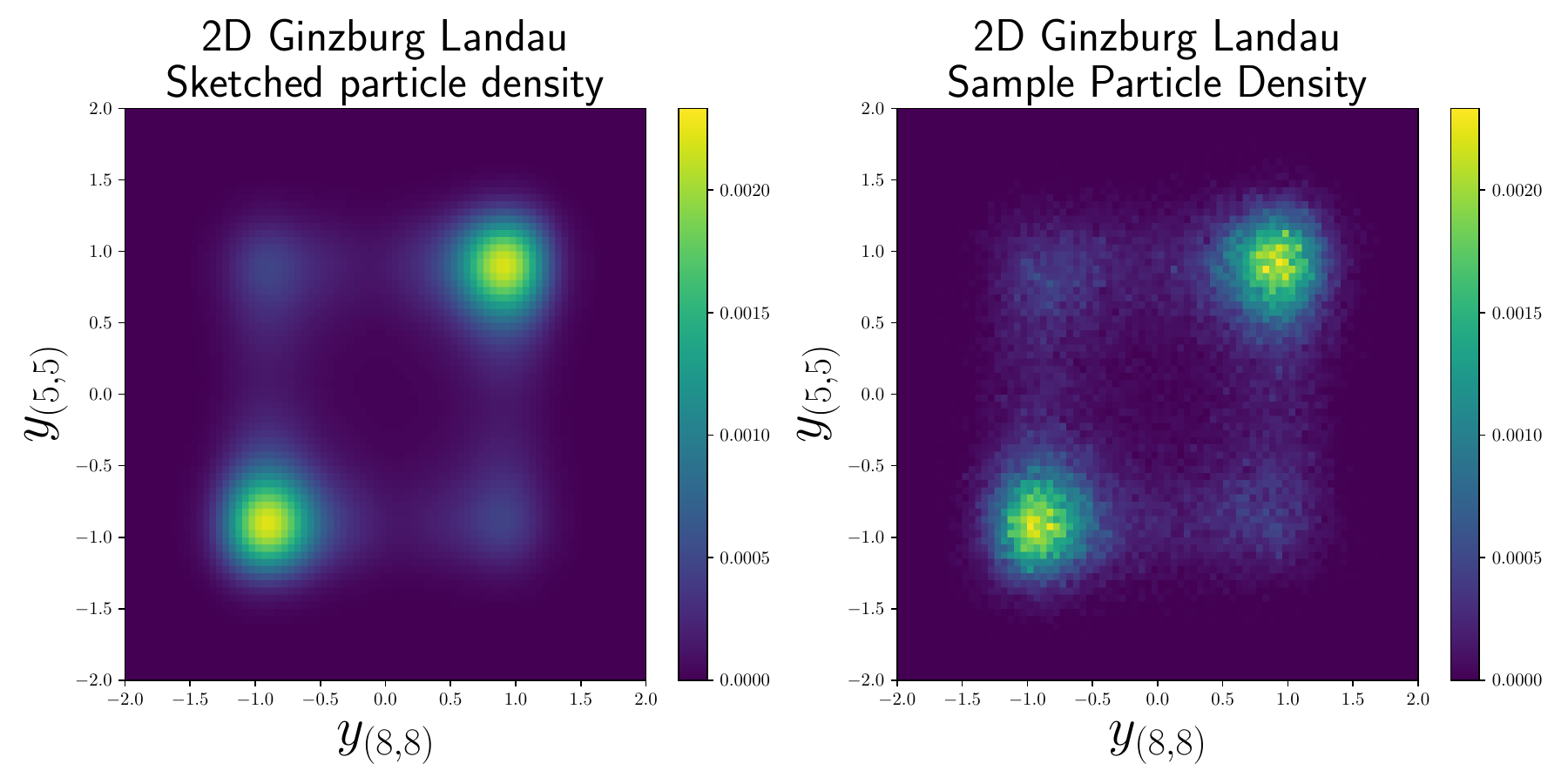}
  \caption{2D Ginzburg-Landau model. The plots of the marginal distribution at \((x_{(5,5)}, x_{(8,8)})\) and \((y_{(5,5)}, y_{(8,8)})\).}
  \label{Fig: 2D marginal}
\end{figure}

In Figure \ref{Fig: 2D marginal}, we plot the marginal distribution of \((x_{(5,5)}, x_{(8,8)})\) and \((y_{(5,5)}, y_{(8,8)})\). In particular, the model can capture the correlation of the \(y\) variable with a relatively small internal bond, which is numerical proof for the hierarchical tensor network ansatz to capture correlation structures when the variables have a 2D lattice structure.

\begin{figure}[!ht]
  \centering
  \includegraphics[width = 0.8\textwidth]{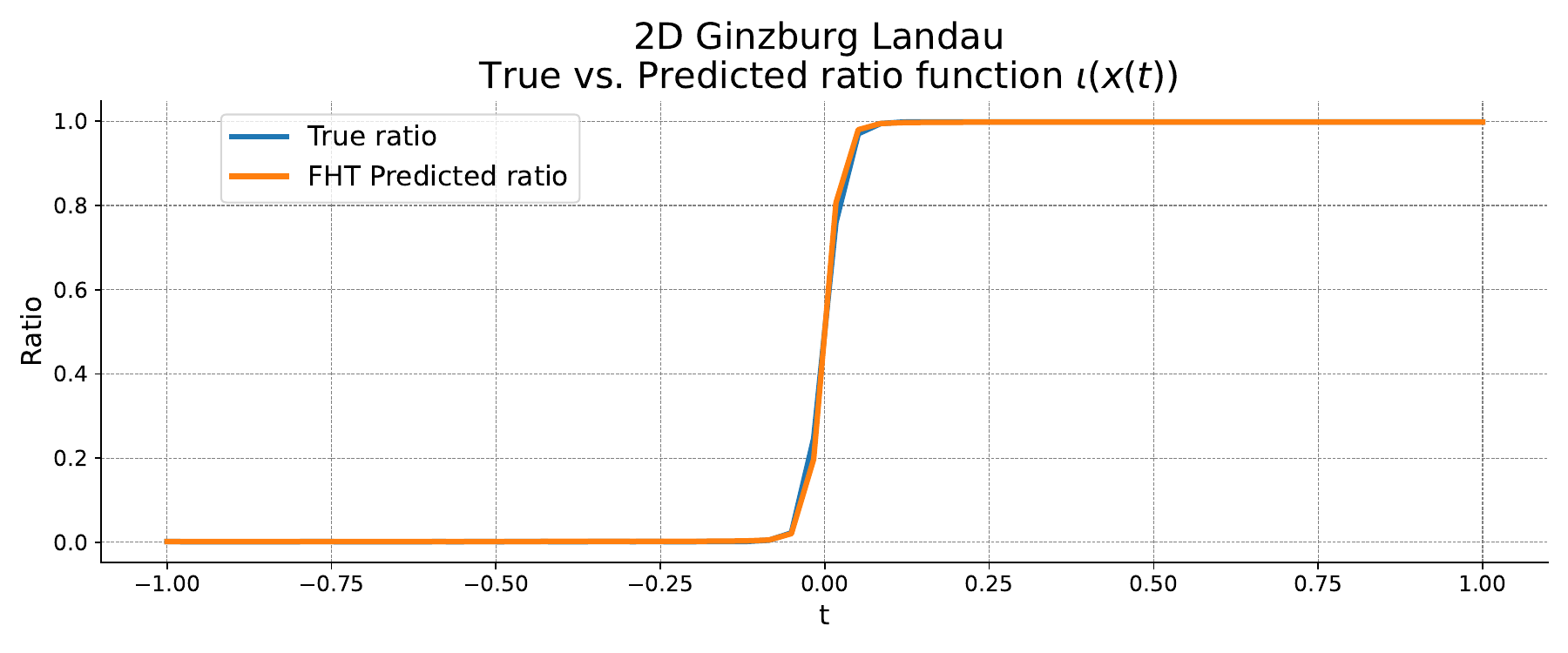}
  \caption{2D Ginzburg Landau model. Plot of \(\iota(x(t))\), i.e. the propensity of \(x(t) = (t,\ldots, t)\) to enter a proximity of \(y_{+} = (1, \ldots, 1)\). One can see that the Markov operator approach can capture the sharp transition at \(t = 0\).}
  \label{Fig: 2D ratio}
\end{figure}

Similar to Section \ref{sec: 1D test}, we consider the KBE solution to terminal conditions \(f_{+}, f_{-}\) as defined in equation \eqref{eqn: terminal condition squared euclidean distance}. Overall, the average relative error for both terminal conditions is less than 8\%. The 2D example has a higher error than the 1D case due to the increased number of variables.
Lastly, in Figure \ref{Fig: 2D ratio}, we repeat the experiment setting in Section \ref{sec: 1D test} and plot \(\iota(x)\) along the direction of \(x(t) = (t, \ldots , t) = t\mathbf{1}_{d}\). One can see that the Markov operator can capture the sharp transition at \(t \in [-0.15, 0.15]\).

\begin{figure}[!ht]
  \centering
  \includegraphics[width = 0.45\textwidth]{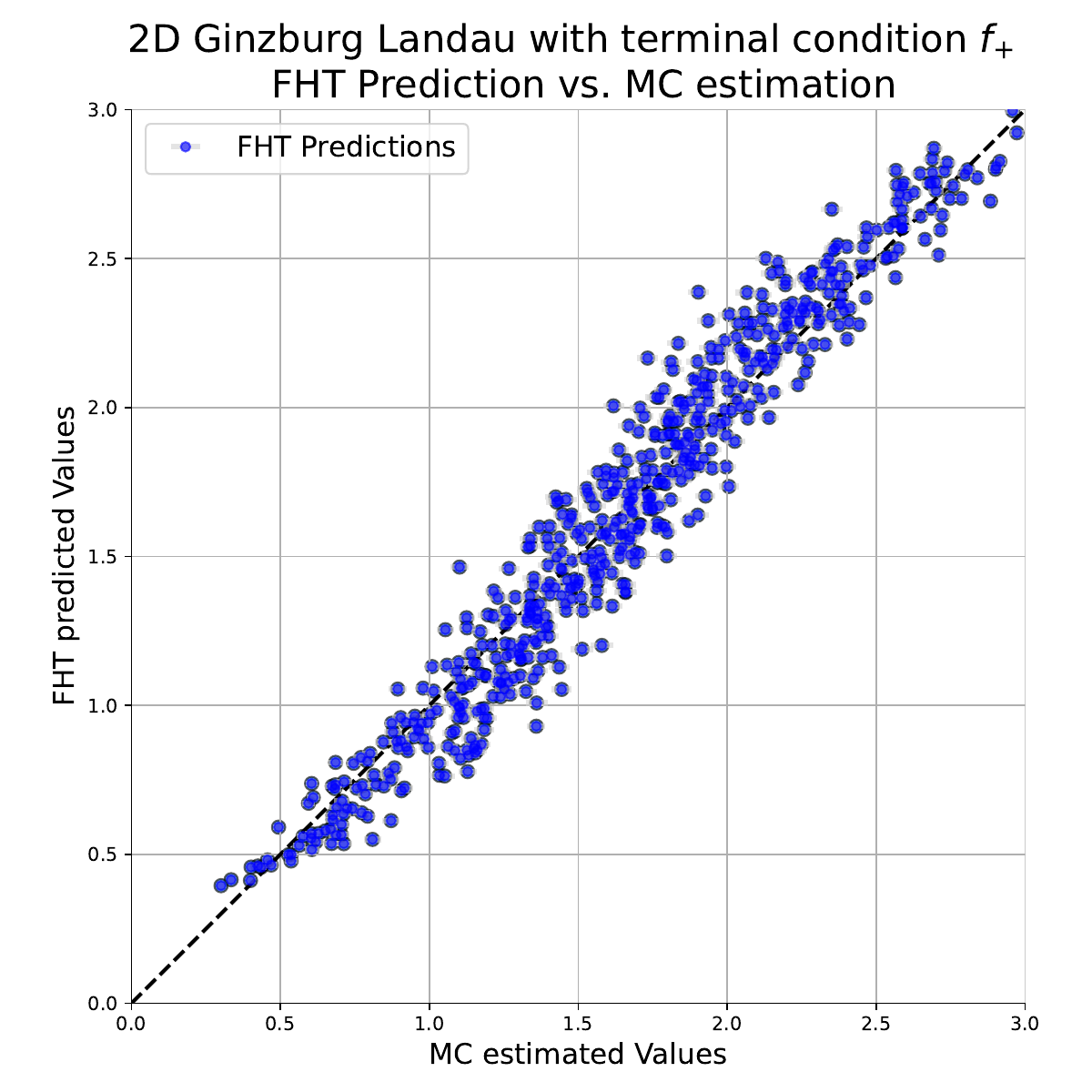}
  \includegraphics[width = 0.45\textwidth]{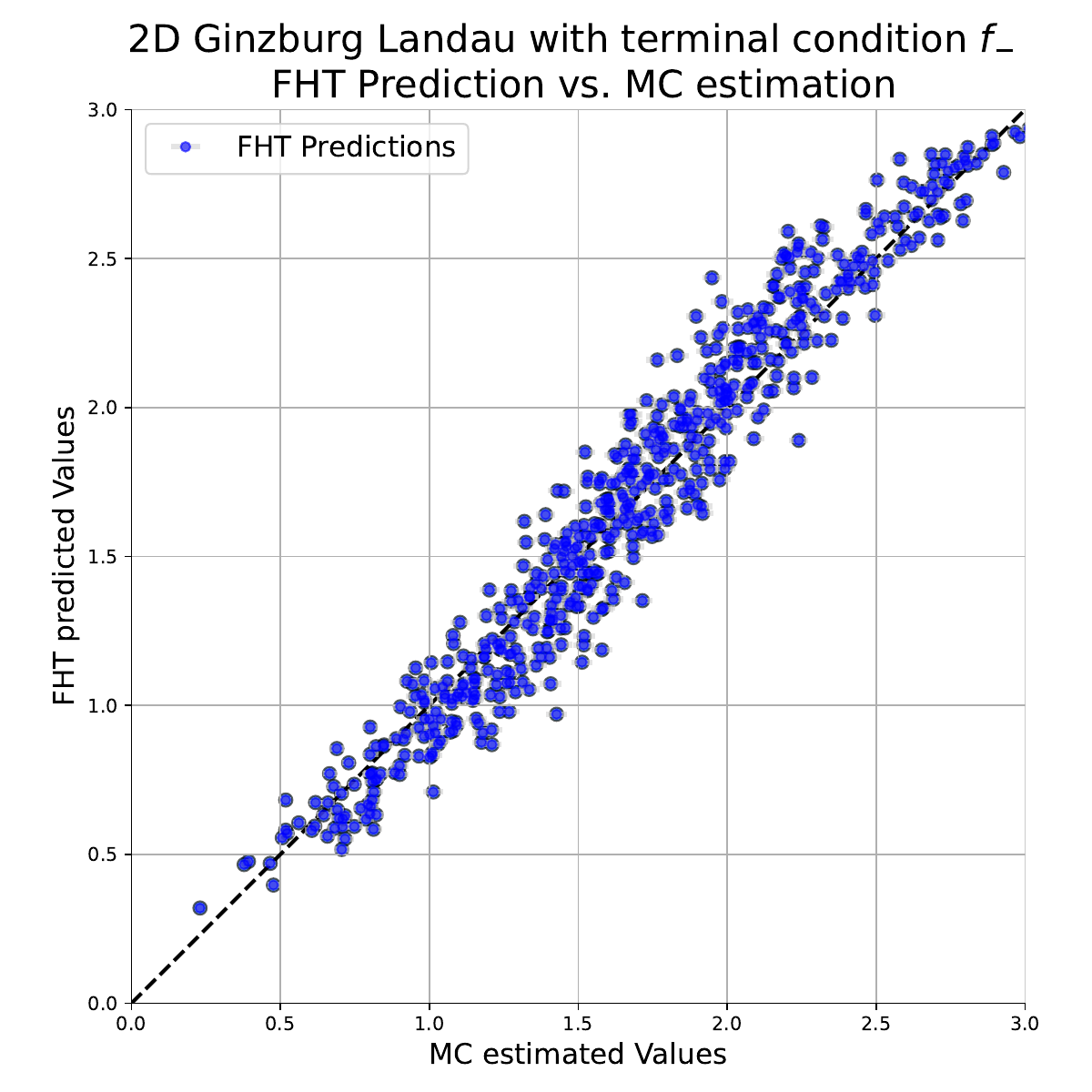}
  \caption{2D Ginzburg Landau model. The plot compares the function evaluation of the FHT-based approximation with reference values obtained from Monte-Carlo estimations, where the error bar indicates the standard deviation from each Monte-Carlo estimation. }
  \label{Fig: 2D error bar}
\end{figure}



\section{Conclusion}
We introduce a novel Markov operator approach that combines particle methods with a functional hierarchical tensor network to solve the high-dimensional Kolmogorov backward equations. The algorithm is applied to the discretized Ginzburg-Landau model in 1D and 2D. This method demonstrates that tensor-network-based operators can be used to tackle high-dimensional equations. An open question is whether one can extend a tensor network approach to solving related nonlinear equations, such as the Hamilton-Jacobi equation.

\bibliographystyle{siam_files/siamplain}
\bibliography{references}

\begin{thebibliography}{10}

\bibitem{bachmayr2016tensor}
{\sc M.~Bachmayr, R.~Schneider, and A.~Uschmajew}, {\em Tensor networks and hierarchical tensors for the solution of high-dimensional partial differential equations}, Foundations of Computational Mathematics, 16 (2016), pp.~1423--1472.

\bibitem{ballani2013black}
{\sc J.~Ballani, L.~Grasedyck, and M.~Kluge}, {\em Black box approximation of tensors in hierarchical tucker format}, Linear algebra and its applications, 438 (2013), pp.~639--657.

\bibitem{beck2021solving}
{\sc C.~Beck, S.~Becker, P.~Grohs, N.~Jaafari, and A.~Jentzen}, {\em Solving the kolmogorov pde by means of deep learning}, Journal of Scientific Computing, 88 (2021), pp.~1--28.

\bibitem{berner2020numerically}
{\sc J.~Berner, M.~Dablander, and P.~Grohs}, {\em Numerically solving parametric families of high-dimensional kolmogorov partial differential equations via deep learning}, Advances in Neural Information Processing Systems, 33 (2020), pp.~16615--16627.

\bibitem{bhattacharyya2019conditional}
{\sc A.~Bhattacharyya, M.~Hanselmann, M.~Fritz, B.~Schiele, and C.-N. Straehle}, {\em Conditional flow variational autoencoders for structured sequence prediction}, arXiv preprint arXiv:1908.09008,  (2019).

\bibitem{brennan1978finite}
{\sc M.~J. Brennan and E.~S. Schwartz}, {\em Finite difference methods and jump processes arising in the pricing of contingent claims: A synthesis}, Journal of Financial and Quantitative Analysis, 13 (1978), pp.~461--474.

\bibitem{cao2018brits}
{\sc W.~Cao, D.~Wang, J.~Li, H.~Zhou, L.~Li, and Y.~Li}, {\em Brits: Bidirectional recurrent imputation for time series}, Advances in neural information processing systems, 31 (2018).

\bibitem{che2018recurrent}
{\sc Z.~Che, S.~Purushotham, K.~Cho, D.~Sontag, and Y.~Liu}, {\em Recurrent neural networks for multivariate time series with missing values}, Scientific reports, 8 (2018), p.~6085.

\bibitem{chen2023committor}
{\sc Y.~Chen, J.~Hoskins, Y.~Khoo, and M.~Lindsey}, {\em Committor functions via tensor networks}, Journal of Computational Physics, 472 (2023), p.~111646.

\bibitem{weinan2004minimum}
{\sc W.~E, W.~Ren, and E.~Vanden-Eijnden}, {\em Minimum action method for the study of rare events}, Communications on pure and applied mathematics, 57 (2004), pp.~637--656.

\bibitem{fan2020solving}
{\sc Y.~Fan and L.~Ying}, {\em Solving electrical impedance tomography with deep learning}, Journal of Computational Physics, 404 (2020), p.~109119.

\bibitem{ginzburg2009theory}
{\sc V.~L. Ginzburg, V.~L. Ginzburg, and L.~Landau}, {\em On the theory of superconductivity}, Springer, 2009.

\bibitem{hackbusch2012tensor}
{\sc W.~Hackbusch}, {\em Tensor spaces and numerical tensor calculus}, vol.~42, Springer, 2012.

\bibitem{hackbusch2009new}
{\sc W.~Hackbusch and S.~K{\"u}hn}, {\em A new scheme for the tensor representation}, Journal of Fourier analysis and applications, 15 (2009), pp.~706--722.

\bibitem{hinton2002training}
{\sc G.~E. Hinton}, {\em Training products of experts by minimizing contrastive divergence}, Neural computation, 14 (2002), pp.~1771--1800.

\bibitem{hoffmann2012ginzburg}
{\sc K.-H. Hoffmann and Q.~Tang}, {\em Ginzburg-Landau phase transition theory and superconductivity}, vol.~134, Birkh{\"a}user, 2012.

\bibitem{hohenberg2015introduction}
{\sc P.~C. Hohenberg and A.~P. Krekhov}, {\em An introduction to the ginzburg--landau theory of phase transitions and nonequilibrium patterns}, Physics Reports, 572 (2015), pp.~1--42.

\bibitem{hur2023generative}
{\sc Y.~Hur, J.~G. Hoskins, M.~Lindsey, E.~M. Stoudenmire, and Y.~Khoo}, {\em Generative modeling via tensor train sketching}, Applied and Computational Harmonic Analysis, 67 (2023), p.~101575.

\bibitem{kovachki2021neural}
{\sc N.~Kovachki, Z.~Li, B.~Liu, K.~Azizzadenesheli, K.~Bhattacharya, A.~Stuart, and A.~Anandkumar}, {\em Neural operator: Learning maps between function spaces}, arXiv preprint arXiv:2108.08481,  (2021).

\bibitem{kumar2006solution}
{\sc P.~Kumar and S.~Narayanan}, {\em Solution of fokker-planck equation by finite element and finite difference methods for nonlinear systems}, Sadhana, 31 (2006), pp.~445--461.

\bibitem{lecun2006tutorial}
{\sc Y.~LeCun, S.~Chopra, R.~Hadsell, M.~Ranzato, and F.~Huang}, {\em A tutorial on energy-based learning}, Predicting structured data, 1 (2006).

\bibitem{lin2011fast}
{\sc L.~Lin, J.~Lu, and L.~Ying}, {\em Fast construction of hierarchical matrix representation from matrix--vector multiplication}, Journal of Computational Physics, 230 (2011), pp.~4071--4087.

\bibitem{lu2020structured}
{\sc Y.~Lu and B.~Huang}, {\em Structured output learning with conditional generative flows}, in Proceedings of the AAAI Conference on Artificial Intelligence, vol.~34(04), 2020, pp.~5005--5012.

\bibitem{mirza2014conditional}
{\sc M.~Mirza and S.~Osindero}, {\em Conditional generative adversarial nets}, arXiv preprint arXiv:1411.1784,  (2014).

\bibitem{peng2023generative}
{\sc Y.~Peng, Y.~Chen, E.~M. Stoudenmire, and Y.~Khoo}, {\em Generative modeling via hierarchical tensor sketching}, arXiv preprint arXiv:2304.05305,  (2023).

\bibitem{rezende2015variational}
{\sc D.~Rezende and S.~Mohamed}, {\em Variational inference with normalizing flows}, in International conference on machine learning, PMLR, 2015, pp.~1530--1538.

\bibitem{richter2021solving}
{\sc L.~Richter, L.~Sallandt, and N.~N{\"u}sken}, {\em Solving high-dimensional parabolic pdes using the tensor train format}, in International Conference on Machine Learning, PMLR, 2021, pp.~8998--9009.

\bibitem{sohn2015learning}
{\sc K.~Sohn, H.~Lee, and X.~Yan}, {\em Learning structured output representation using deep conditional generative models}, Advances in neural information processing systems, 28 (2015).

\bibitem{song2021maximum}
{\sc Y.~Song, C.~Durkan, I.~Murray, and S.~Ermon}, {\em Maximum likelihood training of score-based diffusion models}, Advances in Neural Information Processing Systems, 34 (2021), pp.~1415--1428.

\bibitem{song2019generative}
{\sc Y.~Song and S.~Ermon}, {\em Generative modeling by estimating gradients of the data distribution}, Advances in Neural Information Processing Systems, 32 (2019).

\bibitem{tabak2010density}
{\sc E.~G. Tabak and E.~Vanden-Eijnden}, {\em Density estimation by dual ascent of the log-likelihood}, Communications in Mathematical Sciences, 8 (2010), pp.~217--233.

\bibitem{tang2023solving}
{\sc X.~Tang and L.~Ying}, {\em Solving high-dimensional fokker-planck equation with functional hierarchical tensor}, arXiv preprint arXiv:2312.07455,  (2023).

\bibitem{winkler2019learning}
{\sc C.~Winkler, D.~Worrall, E.~Hoogeboom, and M.~Welling}, {\em Learning likelihoods with conditional normalizing flows}, arXiv preprint arXiv:1912.00042,  (2019).

\bibitem{yoon2018estimating}
{\sc J.~Yoon, W.~R. Zame, and M.~van~der Schaar}, {\em Estimating missing data in temporal data streams using multi-directional recurrent neural networks}, IEEE Transactions on Biomedical Engineering, 66 (2018), pp.~1477--1490.

\bibitem{zhao2007compact}
{\sc J.~Zhao, M.~Davison, and R.~M. Corless}, {\em Compact finite difference method for american option pricing}, Journal of Computational and Applied Mathematics, 206 (2007), pp.~306--321.

\end{thebibliography}

\end{document}